\documentclass[12pt,leqno]{article} 

\usepackage{amsmath} 
\usepackage{amscd} 
\usepackage{amssymb} 
\usepackage{amsfonts}
\usepackage{amsthm}
\usepackage{verbatim} 
\usepackage{pictex}
 \newlength{\fixboxwidth}
\newtheorem{thm}{Theorem}

\newtheorem{rem}{Remark}
\newtheorem{definition}{Definition}
\newtheorem{lemma}{Lemma}
 
\newtheorem{proposition}{Proposition}

\def\e{\varepsilon } 
\def\epsilon{\varepsilon } 
 
\def\rho{\varrho } 
 
\def\phi{\varphi }
\def\a{\alpha }

\def\({\biggl( } 
\def\){\biggr) }

\def\span{{\rm span\,}} 
\def\out{{{\rm out}}}

\def\sq2{{\sqrt{2}}} 
\def\sgn{{\rm sgn }} 
 
\def\A{{\mathcal A}} 
  
\def\lin{{\rm lin}}
\def\non{{\rm non}} 
\def\cont{{\rm cont}}

\newcommand{\Id}{\text{Id}}

\newcount\refnum
\def\refx{\smallskip \global\advance\refnum by 1 {[\the\refnum ] \ }}

\def\a{\alpha}

\def\B{{\cal B}} 
\def\L{{\cal L}} 
\def\N{{\cal N}} 

\def\Id{{\rm id}}
\def\ce{{\mathcal E}}
\def\cl{{\mathcal L}}

\def\C{{\mathcal C}}

\def\D{{\cal D}}
\def\S{{\cal S}} 
\def\Nb{{\mathbb N}} 
\def\R{{\mathbb R}} 
\def\Td{{\mathbb T}^d} 
\def\T{{\mathbb T}} 
\def\Rd{{\mathbb R}^d} 
\def\Cd{{\mathbb C}} 
\def\Z{{\mathbb Z}} 
\def\Zd{{\mathbb Z}^d} 
\def\cf{{\cal F}} 
\def\cfi{{\cal F}^{-1}} 

\def\supp{{\rm supp \, }}

\def\diam{{\rm diam \, }}

\newcommand{\eproof}{\qquad \hfill  \qedsymbol \\}
 
\title{Optimal Approximation of Elliptic Problems by Linear 
and Nonlinear Mappings II} 

\author{Stephan Dahlke\thanks{The work of this author 
has been supported through
the European Union's Human Potential Programme, under contract
HPRN--CT--2002--00285 (HASSIP), and through DFG, 
Grant Da 360/4--2, Da 360/4--3.}, Erich Novak, Winfried Sickel \\ 
\hbox{\small dahlke@mathematik.uni-marburg.de, novak@math.uni-jena.de,
sickel@math.uni-jena.de} }

\begin{document}

\maketitle

\begin{abstract}
We study the optimal approximation of the solution of an operator 
equation
$ 
\A(u) = f 
$ 
by four types of mappings: 
a) linear mappings of rank $n$;
b) $n$-term approximation with respect to a Riesz basis;
c) approximation based on linear information about the
right hand side $f$;
d) continuous mappings. 
We consider 
worst case errors, where $f$ is an element of the unit ball of a 
Sobolev or Besov space
$B^r_q(L_p(\Omega))$ and $\Omega \subset \R^d$ is a bounded
Lipschitz domain; the error is always measured in the $H^s$-norm. 
The respective widths are the linear widths (or approximation numbers),
the nonlinear widths, the Gelfand widths, and the manifold widths.
As a technical tool,  we also study the Bernstein numbers.  
Our main results are the following. 
If $p \ge 2$ then the order of convergence is the same for all 
four classes of approximations. 
In particular, the best linear approximations are of the same 
order as the best nonlinear ones. 
The best 
linear approximation can be quite difficult to realize as a numerical 
algorithm since the optimal Galerkin space usually depends on the 
operator and of the shape of the domain $\Omega$.
For $p<2$ there is a difference, nonlinear approximations 
are better than linear ones. 
However, in this case,  it turns out that linear information 
about the right hand side $f$ is again optimal. 
Our  main theoretical tool is the   best $n$-term approximation 
with respect to an optimal Riesz basis and related nonlinear widths. 
These general results are used to 
study the Poisson equation 
in a  polygonal domain. It turns out 
that best $n$-term wavelet approximation is (almost)
optimal.  
The main results of this paper
are about approximation, not about computation. 
However, we also discuss consequences of the results for the 
numerical complexity of operator equations.
\end{abstract}




\section{Introduction}         


We study the optimal approximation of the solution of an operator 
equation
\begin{equation}     \label{e01} 
\A(u) = f , 
\end{equation} 
where $\A$ is a linear operator 
\begin{equation}     \label{e02} 
\A: H  \to G
\end{equation} 
from a Hilbert space $H$ to another Hilbert space $G$. 
We always assume that $\A$ is boundedly invertible, and so 
\eqref{e01} has a unique solution for any $f \in G$.  
We have in mind the more specific situation of an 
elliptic operator equation which is given as follows. 
Assume that $\Omega \subset \R^d$ is a bounded Lipschitz 
domain and assume that 
\begin{equation}     \label{e03} 
\A : H^s_0 (\Omega) \to H^{-s} (\Omega) 
\end{equation} 
is an isomorphism, where $s > 0$. (For the definition 
of the  Sobolev spaces  $ H^s_0(\Omega)$ and  $H^{-s} (\Omega)$, 
we refer to the
Subsections  \ref{app6}, \ref{app7} and \ref{negative}). 
A standard case (for second order elliptic boundary value problems for 
PDEs) is $s=1$, but also other values of $s$ are of
interest.
Now we put $H=H^s_0(\Omega)$ and $G = H^{-s} (\Omega)$. 
Since $\A$ is boundedly invertible, the inverse mapping 
$S: G \to H$ is well defined. This mapping is sometimes 
called the solution operator---in particular if we want 
to compute the solution $u=S(f)$ from the given 
right-hand side $\A(u)=f$. 

We use linear and (different kinds of) nonlinear mappings $S_n$   for the 
approximation of the solution $u=\A^{-1}(f)$  for $f$ contained in    
$F \subset G$. We consider the worst case error 
\begin{equation}   \label{e04} 
e(S_n, F,H) = 
\sup_{\Vert f \Vert_F \le 1} \Vert \A^{-1}(f)- S_n(f) \Vert_H , 
\end{equation}  
where $F$ is a normed (or quasi-normed) subspace of $G$.
In our main results, $F$ will be a Sobolev or Besov 
space.\footnote{Formally we only deal with Besov spaces. 
Because of the embeddings $B_1^{-s+t}(L_p(\Omega)) \subset 
W_p^{-s+t}(\Omega)  \subset B_\infty^{-s+t}(L_p(\Omega))$, 
which hold for $1 \le p \le \infty$, $t\ge s$, see \cite{T02}, 
our results are valid also 
for Sobolev spaces.} 
Hence we use the following 
commutative diagram
\begin{eqnarray*}
G & \stackrel{\hbox{$S$}}{\longrightarrow} & H \\
I \,  & \nwarrow \qquad  \nearrow & S_F \\
& F. &
\end{eqnarray*}
Here $I: \, F \to G$ denotes the identity and $S_F$ the restriction of 
$S$ to $F$. In the specific case (\ref{e03}) this diagram is given by
\begin{eqnarray*}
H^{-s} (\Omega) & \stackrel{\hbox{$S$}}{\longrightarrow} & H^s_0 (\Omega) \\
I \,  & \nwarrow \qquad  \nearrow & S_t \\
& {B}^{-s+t}_{q} (L_p (\Omega)),  &
\end{eqnarray*}
where ${B}^{-s+t}_{q} (L_p (\Omega))$ denotes a Besov
space compactly embedded into $H^{-s} (\Omega)$, cf. the Appendix 
for a definition, and $S_t$ the restriction of $S$ to 
${B}^{-s+t}_{q} (L_p (\Omega))$.
We are interested in approximations that have the optimal order 
of convergence depending on $n$, where $n$ 
denotes the \emph{degree of freedom}.
In general our results are \emph{constructive in a mathematical sense},
because we can describe optimal approximations $S_n$ in mathematical
terms. This does not mean, however, that these descriptions are 
constructive in a practical sense,  since it might be very difficult 
to convert those descriptions into a practical algorithm. 
We will discuss this more thoroughly in Section~\ref{alg}.
As a consequence, most of our results give  optimal benchmarks 
and can serve for the evaluation of old and new algorithms. 
We study and compare \emph{four kinds of approximation methods};
see Section~\ref{classes} for details. 

\begin{itemize} 

\item
We consider the class $\L_n$ of all continuous linear mappings 
$S_n : F \to H$, 
$$ 
S_n(f) = \sum_{i=1}^n L_i(f) \cdot \tilde  h_i 
$$ 
with arbitrary $\tilde h_i \in H$.  
The worst case error of optimal linear 
mappings is given by the \emph{approximation numbers} 
or \emph{linear widths}
$$ 
e_n^\lin (S , F , H) = \inf_{S_n \in \L_n}  e(S_n, F,H) . 
$$ 

\item
For a given basis $\B$  of $H$ 
we consider the class $\N_n (\B)$ of all (linear or 
nonlinear) mappings of the form 
$$  
S_n(f) = \sum_{k=1}^n c_k \,  h_{i_k} ,      
$$ 
where the $c_k$ and the $i_k$ depend in an arbitrary way on $f$. 
We also allow that the basis $\B$ to be chosen in a nearly arbitrary way. 
Then the \emph{nonlinear widths}  
$e_{n,C}^\non (S,F,H)$ are given by  
$$ 
e_{n,C}^\non (S, F , H) = \inf_{\B \in \B_C} \inf_{S_n \in \N_n(\B)} 
e(S_n, F,H). 
$$ 
Here $\B_C$ denotes a set of Riesz bases for $H$  
where $C$ indicates the stability of the basis. 
These numbers are the main topic of our analysis. 

\item
We also study methods $S_n$ with $S_n = \phi_n \circ N_n$, where 
$N_n: F \to \R^n$ is linear and continuous and 
$\phi_n : \R^n \to H$ is arbitrary. 
This is the class of all (linear or nonlinear) approximations 
$S_n$ that use \emph{linear information of cardinality} $n$
about the right hand side $f$. 
The respective widths are 
$$ 
r_n (S, F, H) := \inf_{S_n} e( S_n, F, H), 
$$
they are closely related to the \emph{Gelfand numbers}. 

\item 
Let $\C_n$ be 
the class of continuous mappings, given by arbitrary 
continuous mappings
$N_n : F \to \R^n$ and $\phi_n : \R^n \to H$. 
Again we define the worst case error of optimal continuous 
mappings by
$$ 
e_n^{\rm cont} (S, F , H) = \inf_{S_n \in \C_n}  e(S_n, F, H ) , 
$$ 
where $S_n = \phi_n \circ N_n$.  
These numbers are called 
\emph{manifold widths} of $S$.

\end{itemize} 

For problems \eqref{e03} with $F=B^r_q(L_p(\Omega))$ 
our main results are the following. 
If $p \ge 2$ then the order of convergence is the same for all 
four classes of approximations. 
In particular, the best linear approximations are of the same 
order as the best nonlinear ones. 
The best 
linear approximation can be quite difficult to realize as a numerical 
algorithm since the optimal Galerkin space usually depends on the 
operator and of the shape of the domain $\Omega$.
For $p<2$ there is an essential difference, nonlinear approximations 
are better than linear ones. 
However, in this case it turns out that linear information 
about the right hand side $f$ is optimal. 
Our  main theoretical tool is     best $n$-term approximation 
with respect to an optimal Riesz basis and related nonlinear widths. 
The main results are about approximation, not about computation. 
However, we also discuss consequences of the results for the 
numerical complexity of operator equations.

\medskip 

The  paper is organized as follows:

\medskip
\noindent 
1. Introduction\\
2. Linear and nonlinear widths\\
2.1 Classes of admissible mappings\\
2.2 Properties of widths and relations between them\\
3. Optimal approximation of elliptic problems\\
3.1 Optimal linear approximation of elliptic problems \\
3.2 Optimal nonlinear approximation of elliptic problems\\
3.3 The Poisson equation\\
3.4 Algorithms and complexity\\
4. Proofs\\ 
4.1 Properties of widths\\
4.2 Widths of embeddings of weighted sequence spaces\\ 
4.3 Widths of embeddings of Besov Spaces\\
4.4 Proofs of Theorems 2, 3, and 5\\ 
5. Appendix - Besov spaces\\
{~}\\
We add a few comments.
The main results of our paper are contained in Section 3.2. 
They are further 
illustrated  for the case of the Poisson equation in Section 3.3. 
A discussion in connection with
\emph{uniform approximation}, 
\emph{adaptive/nonadaptive information}, 
\emph{adaptive numerical schemes},
and
\emph{complexity} 
is contained in Section 3.4. All proofs are contained in Section 4.
Of independent interest are the estimates of the widths of 
embedding operators for Besov spaces, see Section 4.3.
\\
{\bf Notation.} 
We write $a \asymp b$ if there exists  a constant $c>0$ 
(independent of the context dependent 
relevant parameters) such that
\[
c^{-1} \, a \le b \le c \, a \,.
\]
All unimportant constants will be denoted by $c$, 
sometimes with additional indices. 


\section{Linear and Nonlinear Widths} \label{basics} 


Widths represent concepts of optimality. 
In this section we shall discuss several 
variants. Most important for us will be the nonlinear widths $e_n^\non$ and 
the linear widths $e^\lin_n$.
We also study Gelfand and manifold widths and, as a vehicle of 
the proofs, 
Bernstein widths.
  

\subsection{Classes of Admissible Mappings}  \label{classes} 


\subsubsection*{Linear Mappings $S_n$} 

Here we consider the class $\L_n$ of all continuous linear mappings 
$S_n : F \to H$, 
\begin{equation}     \label{e05}  
S_n(f) = \sum_{i=1}^n L_i(f) \,  h_i 
\end{equation} 
where the $L_i: F \to {\mathbb R}$ are linear 
functionals and $h_i $ are elements of 
$H$.
 We consider the worst case error 
\begin{equation}     \label{e06}  
e(S_n, F,H) := 
\sup_{\Vert f \Vert_F \le 1} \Vert \A^{-1}(f)- S_n(f) \Vert_H , 
\end{equation} 
where $F$ is a normed (or quasi-normed) subspace of $G$.
Accordingly, we seek the optimal linear approximation, 
as well as 
the numbers 
\begin{equation}     \label{e07}  
e_n^\lin (S , F , H) = \inf_{S_n \in \L_n}  e(S_n, F,H) , 
\end{equation} 
usually called 
\emph{approximation numbers}  or 
\emph{linear widths}   of $S: F \to H$, 
cf. \cite{Ma90,Pi87,Pi85,Ti90}.     


\subsubsection*{Nonlinear Mappings $S_n$}
 

Let $\B = \{h_1, h_2, \ldots \, \}$ be a subset of $H$.
Then the {\em best $n$-term approximation} of an element $u \in H$ 
with respect to this set $\B$ is defined as
\begin{equation}  \label{e08}  
\sigma_n (u,\B)_H :=
\inf_{i_1, \ldots , i_n}  \inf_{c_1, \ldots \, c_n} \biggl\| u -  
\sum_{k=1}^n c_k \,  h_{i_k} \, \biggr\|_H \, .
\end{equation}  
This subject is widely studied, see the surveys \cite{DV98}
and \cite{Tem03}. Now we continue by looking for an optimal set $\B$ as 
has been done in Kashin \cite{Ka85}, 
Donoho \cite{Do93}, 
Temlyakov \cite{Tem00,Tem02,Tem03} and
DeVore, Petrova, and Temlyakov \cite{DPT03}.  
Temlyakov \cite{Tem03} suggested to consider the quantities 
\[
\inf_{\B \in {\cal D}} 
\sup_{\| u\|_Y \le 1} \, 
\sigma_n (u,\B)_H 
\, , 
\]
where ${\cal D}$ is a subset of the set of all bases of $H$. 
The particular case of ${\cal D}$ being the set of all orthonormal bases
has been discussed in \cite{Tem00,Tem02},
while the set of all unconditional, democratic bases is 
studied in \cite{DPT03}. 
See Remark \ref{VorLimes} for a further discussion. 
In this paper we work with 
Riesz bases, see, e.g.,  Meyer~\cite[page 21]{Me}.

\begin{definition}
Let $H$ be a Hilbert space. Then the sequence $h_1, h_2, \ldots $
of elements of $H$ is called a {\rm Riesz basis}  for $H$ if
there exist positive constants $A$ and $B$ such that, 
for every sequence of scalars
$\alpha_1, \alpha_2, \ldots \, $ with $\alpha_k \not= 0$ for 
only finitely many $k$,  we have
\begin{equation}     \label{e09}  
A \Big(\sum_{k} |\alpha_k|^2 \Big)^{1/2} \le \Big\|
\sum_{k} \alpha_k \, h_k  \Big\|_H \le B 
\Big( \sum_{k} |\alpha_k|^2 \Big)^{1/2} 
\end{equation}
and the vector space of finite sums $\sum \alpha_k \, h_k$
is dense in $H$.  
\end{definition}

\begin{rem} \label{remark1}
The constants $A,B$ reflect the stability of the basis. 
Orthonormal bases are those with
$A=B=1$. Typical examples of Riesz bases are the biorthogonal wavelet 
bases on $\Rd$ or on certain Lipschitz domains, 
cf. {\rm Cohen \cite[Sect. 2.6,~2.12]{C03}}. 
\end{rem}

\noindent
In what follows 
\begin{equation}     \label{e10}  
\B = \{ h_i \mid i \in \Nb \}
\end{equation} 
will always denote a Riesz basis of $H$ with $A$ and $B$ being 
the corresponding 
optimal constants in (\ref{e09}). \\
For a given basis $\B$ 
we consider the class $\N_n (\B)$ of all (linear or 
nonlinear) mappings of the form 
\begin{equation}     \label{e11}  
S_n(f) = \sum_{k=1}^n c_k \,  h_{i_k} ,      
\end{equation} 
where the $c_k$ and the $i_k$ depend in an arbitrary way on $f$. 
By the arbitrariness of $S_n$ one obtains immediately
\begin{equation}   \label{e12}  
\inf_{S_n \in \N_n(\B)} \sup_{\| f\|_F \le 1} \|
\A^{-1} f - S_n (f) \|_H =   
\sup_{\| f\|_F \le 1} \, \sigma_n ( \A^{-1} f,\B)_H \, .  
\end{equation} 
It is natural to assume some common stability 
of the bases under consideration.
For a real number $C \ge 1$ we put
\begin{equation}   \label{e13}
\B_C := \Big\{\B : \, B/A \le C\Big\}.
\end{equation}  
We are ready to define the nonlinear widths
$e_{n,C}^\non (S,F,H)$ by
\begin{equation}     \label{e14}  
e_{n,C}^\non (S, F , H) = \inf_{\B \in \B_C} \inf_{S_n \in \N_n(\B)} 
e(S_n, F,H) . 
\end{equation} 
These numbers are the main topic of our analysis. 
We call them  
the {\em widths of  best $n$-term approximation}
(with respect to the collection $\B_C$ of  Riesz basis of $H$). 

\begin{rem} \label{remark2}
\begin{itemize}
\item[i)]
It should be clear that the class $\N_n (\B)$ contains 
many mappings that are difficult to compute.
In particular, the number $n$ just reflects the dimension of a 
nonlinear manifold and has nothing to do with a computational cost. 
In this paper we also are interested 
in lower bounds, such lower bounds being 
strengthened if we admit a larger 
cass of approximations. 
\item[ii)]
The inequality 
\begin{equation}     \label{e15}  
e_{n,C}^\non   (S, F , H) \le e_n^\lin (S, F , H)  
\end{equation} 
is trivial.
\item[(iii)] Because of the homogeneity of $\sigma_n$, i.e.,
$\sigma_n (\lambda u, \B)_H = |\lambda| \, \sigma_n (u,\B)_H$, 
$\lambda \in {\mathbb R}$, it does not change the 
asymptotic behaviour of $e^\non_n$ if we replace $\sup_{\|f\|_{F} \le 1}$ by
$\sup_{\|f\|_{F} \le c}$ for  $c>0$.
\end{itemize} 
\end{rem}


\subsubsection*{Continuous Mappings $S_n$} 


Linear mappings $S_n$ are of the form $S_n = \phi_n \circ N_n$ where 
both $N_n : F \to \R^n$ and $\phi_n : \R^n \to H$ are linear and 
continuous. If we drop the linearity condition then we obtain 
the class of all continuous mappings $\C_n$, given by arbitrary 
continuous mappings
$N_n : F \to \R^n$ and $\phi_n : \R^n \to H$.
Again we define the worst case error of optimal continuous 
mappings by
\begin{equation}     \label{e16}  
e_n^{\rm cont} (S, F , H) = \inf_{S_n \in \C_n}  e(S_n, F, H ) . 
\end{equation} 
These numbers, or slightly different numbers, were 
studied by different authors, cf. \cite{DHM89,DKLT93,DD96,Ma90}. 
Sometimes these numbers are called 
\emph{manifold widths}  of $S$, see  \cite{DKLT93}, and we will use this 
terminology here. 
The inequality 
\begin{equation}  \label{e17}
e_n^\cont  (S, F , H) \le e_n^\lin (S, F , H)
\end{equation} 
is obvious. 


\subsubsection*{Gelfand Widths and Minimal Radii of Information} 


We can also study methods $S_n$ with $S_n = \phi_n \circ N_n$, where 
$N_n: F \to \R^n$ is linear and continuous and 
$\phi_n : \R^n \to H$ is arbitrary. The respective widths are 
\begin{equation} \label{e18}  
r_n (S, F, H) := \inf_{S_n} e( S_n, F, H). 
\end{equation} 
These numbers are called the $n$-th {\it minimal radii of information}, 
which are closely related to Gelfand widths, see Lemma \ref{bernstein} below.
The {\em $n$-th Gelfand width } of the linear 
operator $S: F \to H$ is given by
\begin{equation}      \label{e19}      
d^n (S, F, H) := \inf_{L_1, \dots , L_n} \, 
\sup \, \Big\{ \Vert Sf \Vert_{H} \, :  \: \Vert f \Vert_F \le 1, \, 
L_i(f) = 0\, , \: i=1, \ldots \, n \Big\} \,  , 
\end{equation} 
where the $L_i : F \to \R$ are continuous linear functionals.


\subsubsection*{Bernstein Widths} 


A well-known tool for deriving lower bounds of widths consists
in the investigation of Bernstein widths,
see \cite{Pi87,Pi85,Ti90}.  

\begin{definition}
The number 
$b_n (S, F , H)$, called the $n$-th {\rm Bernstein width}  of the operator
$S: \, F \to H$,
is the radius of the largest $(n+1)$-dimensional 
ball that is contained in $S(\{ \Vert f \Vert_F \le 1 \})$.
\end{definition}

\begin{rem}\label{Pietsch}
The literature contains several different definitions of Bernstein widths.
For example, Pietsch {\rm \cite{Pi74}} gives the following version.
Let $X_n$ denote subspaces of $F$ of dimension $n$. Then
\[
\widetilde{b}_n (S,F,H) 
:= \sup_{X_n \subset F} \inf_{x \in X_{n}, x\neq 0} \, 
\frac{\| Sx\|_H}{\| x\|_F}\,.
\] 
As long as $S$ is an injective mapping we obviously have
$b_{n} (S,F,H) = \widetilde{b}_{n+1}(S,F,H)$.
\end{rem}


\subsection{Properties of  Widths and Relations Between Them}   
\label{sub2.2}


\begin{lemma}\label{bernstein}
Let $n \in \Nb$ and assume that $F \subset G$ is quasi-normed. \\ 
{\rm (i)} We have $d^n \le r_n \le 2 d^n$ if $F$ is normed and 
$d^n \asymp r_n$ in general. \\ 
{\rm (ii)}
The inequality 
\begin{equation}     \label{e20}  
b_n (S, F, H) \le \min \Big( e_n^\cont (S, F, H), d^n (S,F,H)\Big)
\end{equation} 
holds for all $n$.
\end{lemma}

\begin{rem} 
The inequality $b_n \le e_n^\cont$ is known, compare  
e.g. with {\rm \cite{DHM89}}, and the proof technique 
(via Borsuk's theorem) is often used for the proof 
of similar results. 
\end{rem}

\noindent
The Bernstein widths $b_n$ can also be used to prove lower bounds 
for the $e_{n,C}^\non$. The following inequality has been proved in
\cite{DNS1}.

\begin{lemma}   \label{l2} 
Assume that $F \subset G$ is quasi-normed. 
Then 
\begin{equation}     \label{e21}  
e_{n,C}^\non  (S, F , H) \, \ge \, \frac{1}{2C}  \, 
b_{m} (S, F , H)  
\end{equation} 
holds for all $m \ge 4\, C^2 \, n $.
\end{lemma} 

\noindent
More important for us will be a direct comparison of 
$e^\non_n$ and $e^\cont_n$.
Best $n$-term approximation yields a mapping 
$$
S_n(u) = \sum_{k=1}^n c_k \,  h_{i_k} 
$$
which is in general not continuous. 
However, it is known that certain discontinuous mappings 
can be suitably modified in order to obtain a continuous $n$-term 
approximation with an error which is only slightly worse, 
see, for example, \cite{DKLT93} and \cite{DD00}.  
We prove that, under general assumptions, 
the numbers $e_{n,C}^\non$ can be bounded from below 
by the manifold widths $e_n^\cont$. 

\begin{thm}\label{t1}        
Let $S: G \to H$ be an isomorphism. 
Suppose that the embedding $F\hookrightarrow G$ is compact. Then 
for all $C \ge 1$ and all $n \in \Nb$, we have 
\begin{equation}     \label{e22}     
e_{4n+1}^\cont  (S,F,H) \le 2 \, C \, \| \, S\, \|^2 \, \| \, 
S^{-1}\, \|^2 \, \,  e^\non_{n,C} (S,F,H)\, .
\end{equation}
\end{thm}

\noindent
Finally we collect some further properties of the quantities $e^\cont_n$
and $e^\non_n$.

\begin{lemma}\label{basic} 
{\rm (i)}
Let $m,n \in \Nb$, and let $F$ be a subset of the 
quasi-normed linear space $X$, where  $X$ itself is a subset of the 
quasi-normed linear space  $Y$. Let  $I_j$ denote embedding operators. 
Then 
\begin{equation}\label{Khod}
e_{m+n}^{\rm cont} (I_1,F,Y) \le e_{m}^{\rm cont} (I_2,F,X)\,  
e_{n}^{\rm cont} (I_3,X,Y)
\end{equation}
holds.\\
{\rm (ii)} Let $F$ be a quasi-normed subset of $G$ and let
$I: \, F \to G$ be the embedding. Then
\begin{equation}\label{e24}  
e^\cont_n (I,F,G) \le \| S^{-1}\| \,  
e^\cont_n (S,F,H) \le \| S^{-1}\| \, \| S \| \, 
e^\cont_n (I,F,G)
\end{equation}
and for any $C\ge \| S^{-1}\| \, \| S \| $, we have 
\begin{equation}\label{e25}  
e^\non_{n,C\, \| S^{-1}\| \, \| S \| } (I,F,G) \le \| S^{-1}\| \,  
e^\non_{n,C} (S,F,H) \le \| S^{-1}\| \, \| S \| \, 
e^\non_{n,C/(\| S^{-1}\| \, \| S \| )} (I,F,G) \, . 
\end{equation}
\end{lemma}


\begin{rem} \label{r5}
Let us point out the following which is part of the proof 
of Lemma {\rm \ref{basic}}. Let
$\B= \{h_1, h_2, \ldots \, \}$  be a Riesz basis of $G$. 
Let $S_n$ be an approximation of the identity $I: F \to G$.
Then $S(\B)$ is a Riesz basis of $H$ and $S \circ S_n$ is an 
approximation of $S: F \to H$ satisfying
\begin{equation}  \label{e26}
\Vert f - S_n (f) \Vert_G \le 
\Vert S^{-1} \Vert \cdot \| \, S f - S \circ S_n (f) \, \|_H \le
\Vert S^{-1} \Vert \cdot \| \, S\, \| \cdot  
\| \, f - S_n (f) \, \|_G \, .  
\end{equation}
This makes clear that if $\B$ and  $S_n$ are order optimal
for the triple $I,F,G$, then $S(\B)$ and $S \circ S_n$ are order
optimal for the triple $S,F,H$.
Consequently, instead of looking for good approximations of $S:F \to H$
it will be enough to study approximations of the embedding $I: F \to G$.
\end{rem}

\begin{rem}   \label{r6}  
The assertion in part {\rm (i)}  of the Lemma is essentially 
proved in {\rm \cite{DD96}}  but traced there to Khodulev.
The inequality {\rm (\ref{Khod})}  can be made more 
transparent by means of the diagram
\begin{eqnarray*}
X & \stackrel{\hbox{$I_3$}}{\longrightarrow} & Y \\
I_2 \,  & \nwarrow \qquad  \nearrow & I_1 \\
& F.  &
\end{eqnarray*}
\end{rem}

\begin{rem}\label{r7} 
The
approximation numbers $e^\lin_n$, the
Gelfand widths $d^n$, the manifold widths $e^\cont_n$
and Bernstein widths $b_n$ are particular examples of $s$-numbers
in the sense of {\rm Pietsch \cite{Pi74}}, see 
{\rm \cite{Ma90}} for the manifold widths.
They have several properties in common.
Letting  $s_n$ denote any of the numbers  $e^\lin_n$, $d^n$, $e^\cont_n$ 
and $b_n$ we have
\begin{equation}   \label{mult}
s_n (T_2 \circ T_1 \circ T_0) \le \| \, T_0\, \|\, \, \|\, T_2\, \|
\, s_n (T_1)\, , 
\end{equation} 
where $T_0 \in \cl (E_0,E)$, $T_1 \in \cl (E,F)$, $T_2 \in \cl (F,F_0)$
and $E_0,E,F,F_0$ are arbitrary Banach spaces. For these four types of
$s$-numbers the assertion remains true also for quasi-Banach spaces. 

Another property concerns additivity. For $s_n$ instead of 
$e^\lin_n$ and $d^n$ we have
\begin{equation}   \label{add}
s_{2n} (T_0 + T_1 ) \le c\, \Big(s_n (T_0) + s_n (T_1)\Big)\, , 
\end{equation} 
where $T_0, \, T_1 \in \cl (E,F)$, $E,F$ are arbitrary quasi-Banach spaces,
and $c$ does not depend on $n,T_0,T_1$, cf. {\rm \cite{CS}}.
In case that $F$ is a  Banach space,  one can take $c=1$. 
\end{rem}


\section{Optimal Approximation of Elliptic Problems} 


Let $s, t >0$.
We consider the diagram
\begin{eqnarray*}
H^{-s} (\Omega)  & \stackrel{\hbox{$S$}}{\longrightarrow} & 
{H}^{s}_0 (\Omega)  \\
I \,  & \nwarrow \qquad  \nearrow & S_t \\
& {B}^{-s+t}_{q} (L_{p} (\Omega)),  &
\end{eqnarray*}
where $S_t$ denotes the restriction of $S$ to 
${B}^{-s+t}_{q} (L_{p}(\Omega))$ and $I$ denotes the identity.
We assume \eqref{e03} and we let  $S=\A^{-1}$. 


\subsection{Optimal Linear Approximation of Elliptic Problems} 


\begin{thm}\label{t2} 
Let $\Omega \subset \Rd$ be a bounded Lipschitz domain.
Let $0 < p, q \le \infty$, $s>0$, and 
\begin{equation}    \label{e29}      
t> d \left(  \frac{1}{p} - \frac{1}{2} \right)_+
\, .
\end{equation} 
Then 
\[
e_n^\lin (S,{B}^{-s+t}_{q}  (L_p(\Omega))), {H}^{s}_0 (\Omega) ) 
\asymp \left\{
\begin{array}{lll} n^{-t/d} & \qquad & \mbox{if} \quad 
2 \le p \le \infty\, , \\
n^{- t/d + 1/p -1/2} && \mbox{if}\quad 0< p< 2\, .
\end{array}
\right.
\]
\end{thm}


\begin{rem} \label{remarkacht}
\begin{itemize}
\item[i)] 
The restriction {\rm (\ref{e29})} is necessary and 
sufficient for the compactness of the embedding
$I: B^{-s+t}_q (L_p (\Omega)) \hookrightarrow H^{-s} (\Omega)$, 
cf. the Appendix, 
Proposition {\rm \ref{compact}}.
\item[ii)]
The  proof is constructive. First of all one has to 
determine a linear mapping
$S_n$ that approximates the embedding
$I: B^{-s+t}_q (L_p (\Omega)) \to H^{-s} (\Omega)$ with the optimal order.
How this can be done is described in Remark {\rm \ref{partial}}, 
Subsection {\rm \ref{w34}}.
Finally, the linear mapping $S \circ S_n$ realizes
an in order optimal approximation of $S_t$.
\item[iii)] 
There are hundreds of references dealing with approximation 
numbers of linear operators.
Most useful for us have been the monographs 
{\rm \cite{ET96,Pi87,Pi85,Ti90,Tem93,We96}}, as well as the 
references contained therein. 
\end{itemize}
\end{rem}


\subsection{Optimal Nonlinear Approximation of Elliptic Problems}


To begin with, we consider the manifold and the Gelfand widths. 
There we have a rather final answer.

\begin{thm}\label{t4} 
Let $\Omega \subset \Rd$ be a bounded Lipschitz domain.
Let $0 < p, q \le \infty$, $s>0$, and 
$$
t> d \left(  \frac{1}{p} - \frac{1}{2} \right)_+ \, .
$$
Then 
\[
e_n^\cont (S,{B}^{-s+t}_{q} (L_p (\Omega))), {H}^{s}_0 (\Omega) ) \asymp
n^{-t/d}\, .
\]
If, in addition, $p\ge 1$ (and $t>d/2$ if $1\le p < 2$), then 
\[
d^n (S,{B}^{-s+t}_{q} (L_p (\Omega))), {H}^{s}_0 (\Omega) ) \asymp
n^{-t/d}\, .
\]
\end{thm}

{}From Theorem~\ref{t1} and Theorem~\ref{t4} we conclude that 
the order of $e_{n,C}^\non$ is also at least $n^{-t/d}$. 
For the respective upper bound 
of the nonlinear widths $e_{n,C}^\non$ 
we need a few more restrictions with respect to the domain $\Omega$.
Let $\Omega$ be a bounded Lipschitz domain in $\Rd$ and let $s>0$. 
We assume that for any fixed triple $(t,p,q)$ of parameters
the spaces ${B}^{-s+t}_{q} (L_p (\Omega))$ and 
$H^{-s} (\Omega)$ allow a discretization by one 
common wavelet system $\B^*$, i.e. 
(\ref{eq205})--(\ref{208}) should 
be satisfied with
${B}^{-s+t}_{q} (L_p (\Omega))$ and ${B}^{-s}_{2} (L_2 (\Omega))$,
respectively, 
cf. Appendix \ref{app8}. 
By assumption such a wavelet system belongs to $\B_{C^*}$ for some
$1\le C^* < \infty$.  

\begin{thm}\label{t3} 
Under the above conditions on $\Omega$ and if 
$0 < p, q \le  \infty$, $s >0$,  $t> d(\frac{1}{p} - \frac{1}{2})_+$, 
we have for any $C\ge C^*$
\[
e_{n,C}^\non (S,{B}^{-s+t}_{q} (L_p (\Omega))), {H}^{s}_0 (\Omega) ) 
\asymp  n^{-t/d}\, .
\]
\end{thm}

\begin{rem}   \label{r8} 
Comparing Theorems {\rm \ref{t4}}, {\rm \ref{t3}}   
and Theorem {\rm \ref{t2}} 
there is a clear message.
For $p< 2$ there are
nonlinear approximations that are better in order than any linear 
approximation.
\end{rem}

\begin{rem}   \label{r9}      
The proof of the upper bound in Theorem {\rm \ref{t3}}  is 
\emph{constructive in a theoretical sense} 
that we now describe. 
Given a right-hand side 
$f \in B^{-s+t}_q (L_p(\Omega))$
we have to calculate all wavelet coefficients 
$\langle f,\widetilde{\psi}_{j,\lambda}
\rangle$.
The sequence of these coefficients belongs to the 
space $b^{-s+t}_{p,q}  (\nabla)$, cf. Subsection 
{\rm \ref{sequence}}. With 
\[
a = (a_{j,\lambda})_{j,\lambda}\, , \qquad 
a_{j,\lambda}:= \langle f,\widetilde{\psi}_{j,\lambda} \rangle \, , 
\quad \mbox{for all} \quad j,\lambda\, , 
\]
we find a good approximation 
$S_n(a)$ of $a$ with $n$ components with respect to the norm 
$\| \, \cdot \, | b^s_{2,2} (\nabla)\|$
in Proposition {\rm \ref{proposi*}}.
To get an optimal approximation of the solution $u= Sf$ in 
$\| \, \cdot \, |H^s(\Omega)\|$ we have to apply the solution operator to 
$S_n (a)$. Hence
\begin{equation}\label{e30}  
u_n =  (S \circ S_n) (a) =  \sum_{j=0}^K \sum_{\lambda \in \Lambda_j^*} 
a_{j,\lambda}^* \, S \psi_{j,\lambda}\, ,
\end{equation} 
where $K=K(a,n)$,  with $a_{j,\lambda}^*$ and $\Lambda_j^*$ 
as in Proposition {\rm \ref{proposi*}} (cf. in particular 
{\rm (\ref{Lambda})} 
and {\rm (\ref{lambda3*})}),
represents such a good approximation of $u$.
To calculate $u_n$,  a lot of computations have to be done.
The coefficients $a_{j,\lambda}^*$ are the largest in a weighted sense
(the weight depends on $n$ and $j$, cf. the proof of Proposition 
{\rm \ref{proposi*}} 
for explicit formulas). Having these coefficients at hand one has 
finally to solve all the equations
\begin{equation}   \label{e31}
\A u_{j, \lambda}  = \psi_{j,\lambda}\, , 
\qquad 0 \le j \le K, \quad \lambda \in \Lambda_j^* 
\end{equation}
to obtain $u_{j, \lambda} = S \psi_{j, \lambda}$. 
The number of equations is $O(n)$.

In this way we obtain a nonlinear approximation with respect to the
Riesz basis given by the $S \psi_{j, \lambda}$. Observe 
that this Riesz basis depends on the operator equation. 
It would be much better to use a known Riesz basis, 
such as a wavelet basis, that does not depend on $\A$. 
See Theorem~{\rm \ref{t5}} for a step into that direction. 
\end{rem}

\begin{rem}   \label{r10} 
At least if $\Omega$ is a cube, 
all required properties are known to be satisfied
if in addition  $1<p,q<\infty$. The latter restriction allows to use 
duality arguments, cf. Proposition {\rm \ref{duality}} in Appendix
{\rm \ref{app7}}. 
 There also exist
results for domains with piecewise analytic boundary such as polygonal or\
polyhedral domains. One natural way as, e.g., outlined in {\rm \cite{CTU}} and
{\rm \cite{DS99}}, is to decompose the domain into a  disjoint union 
of parametric
images of reference cubes. Then, one constructs wavelet bases on the reference
cubes and glues everything together in a judicious fashion. However, due to the
glueing procedure, only Sobolev spaces $H^s$ with smoothness $s <3/2$ can be
characterized. This bottleneck can be circumvented by the approach in
{\rm \cite{DS99a}}.  
There,  a much more tricky domain decomposition method involving
certain projection and extension operators is used. By proceeding in 
this way, norm
equivalences for all spaces $B^t_q(L_p(\Omega))$ can be derived, at least for
the case $p>1$, see {\rm \cite[Theorem 3.4.3]{DS99a}}.  However, the 
authors also mention that their results can be generalized to the case $p<1$, 
see {\rm \cite[Remark 3.1.2]{DS99a}}.

Sobolev and Besov spaces on compact $C^\infty$-manifolds were already 
characterized via spline bases and sequence spaces by Ciesielski and 
Figiel {\rm \cite{CF}}. In that paper also the isomorphism between 
function spaces and sequence spaces is used to obtain results 
for various $s$-numbers. 
\end{rem}

\begin{rem} \label{r11}
Comparing Theorems {\rm \ref{t4}} and {\rm \ref{t3}}
we see that the numbers 
$e^\non_{n,C}$, $e_n^\cont$, and $d^n$ have the same asymptotic behaviour, at
least for $p>1$.
Using the relation $d^n \asymp r_n$, see Lemma~{\rm \ref{bernstein}}, 
we actually can get the optimal order $n^{-t/d}$ with 
an approximation of the form
\begin{equation} \label{e32}  
f \mapsto S \circ \phi_n \circ N_n (f) \, , 
\end{equation} 
where 
$$
N_n : B^{-s+t}_{q} (L_p (\Omega)) \to \R^n 
$$
is linear
(this mapping gives the \emph{information} 
that is used about the right hand side), and 
$$
\phi_n :  \R^n \to H^{-s} (\Omega) 
$$
is nonlinear. 
Note that neither $N_n$ nor $\phi_n$ depend on $S$.
The mapping $\phi_n \circ N_n$ gives a good 
approximation of the embedding from 
$B^{-s+t}_{q} (L_p (\Omega))$ to $H^{-s}$. 
\end{rem} 

\begin{rem}
There is a further little difference between linear 
and nonlinear approximation.
Let us consider the limiting case 
$t= d(1/p-1/2)$,  where  $0 < p < 2$. Then the embedding
$B^{-s+t}_p (L_p (\Omega)) \hookrightarrow H^{-s} (\Omega)$ 
is continuous, not compact.
As a consequence 
\[
e_n^\lin (S,{B}^{-s+t}_{p}  (L_p(\Omega))), {H}^{s}_0 (\Omega) ) 
\not \to 0 \qquad 
\mbox{if}\quad n \to \infty\, , 
\]
but
\[
e_n^\non (S,{B}^{-s+t}_{p}  (L_p(\Omega))), {H}^{s}_0 (\Omega) )  
\to 0 \qquad 
\mbox{if}\quad n \to \infty\, ,
\]
cf. Remark {\rm \ref{Limes}}.
\end{rem}

 
\subsection{The Poisson Equation} 


The next step  is to discuss the specific case of
the Poisson equation on  a Lipschitz domain $\Omega$ contained in $\R^2$: 
\begin{eqnarray}
-\triangle u &=&f \quad  \mbox{in}\quad \Omega \label{Poisson}\\
           u&=&0 \quad \mbox{on} \quad \partial \Omega. \nonumber
\end{eqnarray}

\noindent
As usual, we study  (\ref{Poisson}) in the weak formulation. Then, it can be
shown that  the operator ${\cal A}= \triangle:~H_0^1\longrightarrow H^{-1}$ is
boundedly invertible, see, e.g., \cite{H92}  for details. Hence
Theorems~\ref{t2} and~\ref{t4} apply with $s=1$; 
for the upper bound of Theorem~\ref{t3} we need 
some restrictions with respect to $\Omega$. 
For the proof of Theorem~\ref{t3} we used the 
Riesz basis 
$S \psi_{j, \lambda}$,  which depends on $\A$. 
Now we want to approximate the solution $u$ 
by wavelets. 
\\
We shall restrict ourselves to the case that 
$\Omega$  is   a simply connected polygonal 
domain. The segments of $\partial \Omega$ are
denoted
by $\overline{\Gamma}_1, \dots , 
\overline{\Gamma}_N$, where each $\Gamma_l$ is open and 
the segments are 
numbered in positive  orientation. Furthermore, $\Upsilon_l$
denotes the endpoint of $\Gamma_l$ and $\omega_l$ denotes 
the measure of the interior  angle at $\Upsilon_l$. 
Appropriate
wavelet systems can be  constructed for such a domain, see Remark \ref{r10}.
Then we obtain the following.

\begin{thm}    \label{t5}
Let $\Omega $ be a polygonal domain 
in $\R^2$. Let  $1 < p \le  2$  and let $ k\geq 1$  
be a nonnegative integer such that
$$ 
\frac{m\pi}{\omega_l} \not=k+1-\frac{2}{p} 
\qquad \mbox{for all} \quad m \in {\mathbb N}, ~l=1, \ldots, N.$$
Then for an appropriate wavelet system $\B^*$,  
the best $n$-term approximation of
problem {\rm (\ref{Poisson})}   yields 
\begin{equation}    \label{e34}   
\sup_{\Vert f \mid B^{k-1}_p(L_p(\Omega)) \Vert  \le 1}
\sigma_n (u, \B^* )    \leq c_\e \, n^{-k/2 + \epsilon}
\end{equation}
where $\varepsilon >0$ and $c_\e$ do not depend on $n$. 
\end{thm}


\begin{rem}  \label{r12} 
This approximation differs greatly from the one 
described in Remark~{\rm \ref{r9}}.
Here we can work with one given wavelet system to approximate 
the solution $u$. We are not forced 
to work with the solutions of the system {\rm (\ref{e31})}. A more detailed
discussion of these relationships, 
including  possible numerical realizations 
of wavelet methods,  will follow in Section {\rm \ref{alg}}.  
\end{rem}

 
\subsection{Algorithms and Complexity}   \label{alg} 


So far,  we have studied the error $e(S_n, F, H)$ of approximations $S_n$. 
We compared the error of nonlinear $S_n$ and linear $S_n$ 
and proved results on the optimal rate of convergence. 
We assume that \eqref{e01} is a given  fixed operator equation
and hence, in the case of \eqref{e03}, also $\Omega$ is fixed. 

In this section we briefly discuss algorithms
and their complexity,  and for simplicity 
we still assume that the operator equation \eqref{e03} 
is given and fixed. Observe that in practice it is important to 
construct also algorithms for more general problems: 
We want to input information about $\Omega$ and $\A$ and the right 
hand side $f$,  and  we want to obtain an $\e$-approximation of the solution
$u$. In our more restricted case we \emph{only} have to input 
information concerning the right hand side $f$ because $\Omega$ and
$\A$ are fixed. 

As is usual in numerical analysis, we use the real 
number model of computation
(see \cite{N95} for the details and \cite{NW99} and \cite{NW00} 
for further comments).
Any algorithm computes and/or uses some information 
(consisting in finitely many numbers) describing the 
right hand side $f$ of \eqref{e03}. 
There are different ways how an algorithm may use information 
concerning $f$, we describe two of them in turn. 

\begin{enumerate} 

\item
The information used about $f$ 
is very explicit if $S_n$ 
is linear \eqref{e05}: 
Then the algorithm uses 
$L_1(f), \dots , L_n(f)$ and we assume that we have an oracle 
(or subroutine) for the $L_i(f)$. 
In practical applications the computation of a functional 
$L_i(f)$ can be very easy or very difficult or anything between. 
One often assumes that the cost of obtaining a value $L_i(f)$ is 
$c$ where $c >0$ is small or large, depending on the circumstances. 

As in \eqref{e11}, we can imagine $S_n$ as the input-output mapping 
of a numerical algorithm: on input $f \in F$ we obtain the output 
$
S_n (f) = u_n = \sum_{k=1}^n c_k \, h_{i_k} 
$. 
More formally we should say that the output is 
\begin{equation}    \label{e37}  
\out (f) = ( i_1, c_1, i_2, c_2, \dots , i_n , c_n ) 
\end{equation} 
but we identify $\out (f)$ with $u_n$. 
Of course we cannot consider arbitrary mappings 
$S_n$ of the form \eqref{e11}  
as the input-output mapping of an algorithm, 
since not all such $S_n$ are computable. 

%
%
%
%
%

We still assume that we only have an oracle for the computation 
of linear functionals $L_i(f)$. 
Then it is not so clear what the information cost of \eqref{e11}  is, since 
\eqref{e11}  only describes the (desired) output of an algorithm, it is 
not an algorithm by itself. We need an algorithm that uses information 
$L_1 (f) , \dots , L_N (f)$, where $N$ might be bigger than $n$, to produce
the $i_k$ and the $c_k$ of $\out (f)$. The information cost of such a 
procedure would be $c N$. 

\item
One also can assume that a good approximation $f_n$  
can  easily be precomputed with negligible cost. 
Hence 
the algorithm starts with an approximation 
\begin{equation} \label{e38} 
f_n = \sum_{k=1}^n c_k \, g_{i_k} ,
\end{equation} 
such as a best $n$-term approximation (or a greedy 
approximation) of $f$ with respect 
to a basis $\{ g_i, \, : \,  i \in \Nb \}$. 

\end{enumerate} 

This is a good place for a short remark about adaption. The use of 
\emph{adaptive methods} is quite widespread but we want to stress that the 
notion of adaptive methods is not uniformly used in the literature. 
Some confusion is almost unavoidable if such different notions are 
mixed. To avoid such confusion,  we do not use the notion of an 
``adaptive method''. Instead we speak 
first about \emph{adaptive (or nonadaptive) 
information} and then about 
\emph{adaptive numerical schemes}.

\begin{itemize}   

\item
\emph{Nonadaptive information}:
The algorithm uses certain functionals 
$L_1, L_2, \dots , L_n$ and for each input $f \in F$ the algorithm needs 
$L_1(f), L_2(f), \dots , L_n (f)$. Hence the functionals $L_i$ 
do not depend on $f$. 
In this case we say that the algorithm uses \emph{nonadaptive information}. 

\item
\emph{Adaptive information}: 
The algorithm uses $L_1(f)$ and, depending on this number, the next 
functional $L_2$ is chosen. In general, the chosen functional 
$L_k$ may depend on the values $L_1(f), \dots , L_{k-1}(f)$ 
that are already known to the algorithm. 
Observe that $L_k$ cannot depend in an arbitrary way on $f$ since 
the algorithm can only use the known information about $f$. 
In this case we say that the algorithm uses 
\emph{adaptive information}. 

\end{itemize}   

We give an example. 
Assume that a certain $S_n$ of the form \eqref{e11}  can be realized
in such a way that we first compute $L_1(f), \dots , L_N (f)$,
where the $L_i$ do not depend on $f \in F$. In the latter parts 
of the algorithm we only use the $L_i(f)$ for the $n$ largest values of 
$|L_i(f)|$, together with the corresponding values of $i$, to compute the 
output $\out (f)$. Such an algorithm uses nonadaptive information
(of cardinality $N$),
the information cost is $c N$. 

There is a large stream of results, 
giving conditions under which adaptive
information is superior (or not superior) compared to nonadaptive
information; we mention the pioneering paper by
Bakhvalov \cite{B71}, the results on operator equations 
by Gal and Micchelli \cite{GM80} and
by Traub and Wo\'zniakowski \cite{TW80},  
and the survey \cite{N96}. 
For example, it is known that adaptive information does not help 
(up to a factor of 2) for linear operator equations
and the worst case 
error with respect to the unit ball of a normed space $F$. If $F$ is only 
quasi-normed then the proofs must be modified, with a possible change of the
constant 2. Nevertheless nonadaptive information is almost as good as 
adaptive information.

How much information is needed about the right hand side $f \in F$
in order that we can solve the equation \eqref{e01} with an error $\e$?
This question is answered by 
the minimal radii of information $r_n(S,F,H)$ 
(or the closely related Gelfand numbers). These numbers are a good 
measure for the \emph{information complexity} of the operator equation.
In contrast, the \emph{output complexity} of the problem is measured 
by the nonlinear widths $e_{n,C}^\non (S, F, H)$. These numbers measure 
the cost of just outputting the approximation 
(with respect to an optimal basis $\B \in \B_C$). 
It is quite remarkable that, under general conditions, 
we obtain the same order
$$
r_n(S,F,H)  \asymp d^n(S,F,H)  \asymp e_{n,C}^\non(S,F,H)  \asymp n^{-t/d},
$$
see Theorem~\ref{t4} and Theorem~\ref{t3}. 

Now we discuss 
\emph{adaptive numerical schemes}
for the numerical treatment of  
elliptic partial differential equations.
Usually, these operator equations are solved by a Galerkin scheme, i.e.,
one defines an increasing sequence of finite 
dimensional approximation spaces $G_{\Lambda_l}
:=\mbox{span}\{\eta_{\mu}: \mu \in \Lambda_l\},$ where    
$G_{\Lambda_l} \subset G_{\Lambda_{l+1}}$,   
and projects  the problem onto these spaces, i.e.,
$$
\langle {\cal A}u_{\Lambda_l},v\rangle= \langle f ,v\rangle   
\quad \mbox{for all} \quad v \in G_{\Lambda_l} . 
$$

To compute the actual Galerkin approximation, 
one has to solve a linear system
$${\bf A}_{\Lambda_l} {\bf c}_{\Lambda_l}= {\bf f}_{\Lambda_l}, 
\qquad  \qquad{\bf A}_{\Lambda_l}=
(\langle A\eta_{\mu'}, \eta_{\mu}\rangle)_{\mu,\mu' 
\in \Lambda_l}, \qquad ({\bf f}_{\Lambda})_{\mu}=\langle f,\eta_{\mu}
\rangle,~ \mu \in \Lambda_l.$$
Then the question arises how to choose the approximation spaces  
in  a suitable way, since doing that in a somewhat clumsy  fashion
would yield huge linear systems and a very unefficient scheme. 
One natural way would be to use an updating strategy,
i.e., one starts with a small set $\Lambda_0$, 
tries to estimate the (local) error, and only in regions where
the error is large the index set is {\it refined}, i.e., 
further basis functions are added. Such an updating strategy is
usually called an {\it adaptive numerical scheme} 
and it is characterized by the following facts: the sequence 
of approximation spaces is not a priori fixed but depends  
on the {\it unknown} solution $u$ of the operator equation,
and the whole scheme should be self-regulating, i.e., 
it should work without a priori information 
on the solution.   In principle,
such an adaptive scheme consists of the following three steps:
$$ 
\begin{array} {ccccc}
\mbox{solve}&-&\mbox{estimate}&-&\mbox{refine}\\[2mm]
{\bf A}_{\Lambda_l}{\bf c}_{\Lambda_l}={\bf f}_{\Lambda_l} 
&& \|u-u_{\Lambda_l}\|=?&& \mbox{add  functions}\\
&& \mbox{a posteriori} && \mbox{if necessary.}\\
&&\mbox{error estimator} &&
\end{array}
$$
Note that the second step is highly nontrivial since the exact 
solution $u$ is unknown, so that
clever a posteriori error estimators are needed. 
These error estimators should be local,  since
we want to refine (i.e. add basis functions) only in regions 
where the local error is large.
Then another
challenging task is to show  that the refinement 
strategy leads to a convergent scheme and to estimate
its order of convergence, if possible. 

Recent developments
indicate the promising potential  
of adaptive numerical schemes, see, e.g.,  \cite{BR2,BW,Ra1,BEK,Do,Steve,Ve2}
for finite element methods. 
However, to further explain the ideas and to make 
comparisons as simple as possible, we shall restrict ourselves 
to adaptive schemes based on wavelets.  
For simplicity, we shall  mainly discuss  the approach in \cite{DDHS};
for more sophisticated versions the reader 
is referred to \cite{CDD,CDD02,CDD3, DDU}.  
The first step clearly  must be the development 
of an a posteriori error estimator. Using
the fact that ${\cal A}$ is boundedly invertible 
and the  usual norm equivalences, compare with (\ref{208}),  we obtain

\begin{eqnarray}
\|u-u_{\Lambda}\|_{H^s}&\asymp  &
\|{\cal A}(u-u_{\Lambda})\|_{H^{{-s}}}  \label{e39} \\
&\asymp  & \|f-{\cal A}(u_{\Lambda})\|_{H^{-s}} \nonumber\\
& \asymp  & \|r_{\Lambda}\|_{H^{-s}} \nonumber \\
&\asymp  & \biggl(  \sum_{{(j,\lambda) }\in J 
\backslash \Lambda} 2^{-2{s}{j}}|\langle r_{\Lambda}, 
\psi_{{j,\lambda}}\rangle|^2 \biggr)^{1/2} \nonumber \\
&=& \biggl(\sum_{{(j,\lambda)}\in J \backslash \Lambda}  
{\delta_{{j,\lambda}}^2}\biggr)^{1/2}, \nonumber
\end{eqnarray}
where the {\it residual weights} $\delta_{{j,\lambda}} $ can be computed as
$$
\delta_{{j,\lambda}}= 2^{-s{j}} \biggl| f_{{j,\lambda}}-\sum_{{
(j',\lambda')} \in \Lambda} \langle {\cal A} \psi_{{j',\lambda'}},
\psi_{{j,\lambda}}\rangle u_{{ j',\lambda'}} \biggr|  
\quad \hbox{with} \quad 
f_{j,\lambda}=\langle f, \psi_{{j,\lambda}}
\rangle.
$$ 
{}From (\ref{e39}),  we observe that the 
sum of the residual weights gives rise
to an efficient and reliable a posteriori error estimator. 
Each residual weight $\delta_{j,\lambda}$ can be 
interpreted as a local error indicator, so that the following 
natural refinement strategy suggests itself: 
Add wavelets in regions where the residual weights are large; 
that is, try to catch the bulk of the residual
expansion in (\ref{e39}). 
Indeed, it can be shown that this strategy produces a convergent adaptive
scheme, in principle. However, we are faced with a serious problem: 
the index set $J$ will not have 
finite cardinality, so that  neither the  error estimator   
nor the adaptive refinement strategy can  be implemented.
Nevertheless, there exist implementable variants, 
see again \cite{CDD,DDHS} for details. 
%
%
We  start  with  the set
$$J_{{j,\lambda}, {\varepsilon}}:~\{ {(j',\lambda')}|~|\langle {\cal A}
\psi_{{j',\lambda'}}, \psi_{{j,\lambda}}\rangle|~{\varepsilon}
\mbox{-significant}\}$$
and define
$$a_{{j,\lambda}}(\Lambda,  {\varepsilon}):=2^{{-s}{j}}|
\sum_{{(j',\lambda')} \in \Lambda \cap
 J_{{j,\lambda},{\varepsilon}}}\langle {\cal A} \psi_{{j',\lambda'}},
\psi_{{j,\lambda}}\rangle u_{{j',\lambda'}}|.$$  
(The expresssion  `$\varepsilon$-significant' can be made 
precise by using the locality and
the cancellation properties of a wavelet basis).  
By employing the $a_{{j,\lambda}}(\Lambda,  {\varepsilon})$
we obtain another error erstimator:
$$
\|u-u_{\Lambda}\|_{H^s} 
\le 
c \cdot \biggl( \biggl(\sum_{{(j, \lambda)} \in J \backslash \Lambda}
a_{{j,\lambda}}^2\biggr)^{1/2} \!\!\!\!\!\!\!\!+ 
{\varepsilon}\|f\|_{H^{{-s}}} +
\!\!\inf_{v \in \tilde{V}_{\Lambda}}\!\|F-v\|_{H^{{-s}}} \biggr) . $$
Here $\tilde{V}_{\Lambda}$ denotes the approximation 
space spanned by the dual wavelets corresponding to
$\Lambda$, see Section \ref{app3}  for details. 
Now, playing the same game for the $a_{{j,\lambda}}(\Lambda,  
{\varepsilon})$ instead of the
$\delta_{j,\lambda}$, we end up with a convergent and 
implementable adaptive strategy. 
To this end, the starting index set $\Lambda$  has to be 
determined such that  $\inf_{v \in \tilde{V}_{\Lambda}}
\|f-v\|_{H^{{-s}}} \le  
c \cdot {\mbox{eps}}$ and  ${ \varepsilon}(f,\mbox{eps}, \theta)$ 
has  to be computed. Then, there exists  
%
%
a  constant ${\kappa} \in (0,1)$ such that  whenever 
$\tilde{\Lambda} \subset J,~ 
\Lambda \subset \tilde{\Lambda}$ is chosen so that  
\begin{equation} \label{e40}   \left(\sum_{{(j,\lambda)} 
\in \tilde{\Lambda} \backslash \Lambda} a_{j,\lambda}(\Lambda, 
{\varepsilon})^2 \right)^{1/2} \geq (1-{\theta})
\left( \sum_{{(j,\lambda)} \in J\backslash \Lambda} a_{{j, \lambda}}
(\Lambda, {\varepsilon})^2 \right)^{1/2}\end{equation}
 either
\begin{equation} \label{energyred} 
\|u-u_{\tilde{\Lambda}}\| \leq \kappa\| u - u_{\Lambda}\|, 
\qquad {\kappa} \in (0,1)
\end{equation}
or
\begin{equation} \label{e41}   \left(\sum_{(j,\lambda) \in J 
\backslash \Lambda} a_{{j,\lambda}}(\Lambda, { \varepsilon})^2 \right)^{1/2}
\leq { \mbox{eps}} \end{equation}
which implies that \begin{equation} \label{e42}    \|u-u_{\Lambda}\| \le   
{\mbox{eps}} \cdot c  .
\end{equation} 
For the proof and further details, the reader is again refered to
\cite{DDHS}.
\begin{rem}  \label{r13} {
\begin{itemize}
\item[i)]  In order to avoid unnecessary technical  
and notational difficulties,
we have not presented the explicit form of the function 
$\varepsilon( f, \mbox{eps}, \theta)$. It depends in a complicated, but
nevertheless computable way on the  final accuracy $\mbox{eps}$, the control
parameter $\theta$ , the $H^{-s}$--norm of the right--hand side $f$, and on
the stability and ellipticity constants of the problem. For details, we refer
again to {\rm \cite{DDHS}}.
\item[ii)]  The norm $\|\cdot\|$ in {\rm (\ref{energyred})} and 
{\rm (\ref{e42})} clearly denotes the energy norm $\|v\|:=\langle {\cal A}v, v
\rangle$,  which is equivalent to the Sobolev norm $H^s$, see again {\rm
\cite{H92}}  for details. 
\item[iii)]  Eqs. {\rm (\ref{energyred}), 
(\ref{e41})} and {\rm (\ref{e42})}  obviously imply
that the adaptive strategy  in {\rm (\ref{e40})} converges. Indeed,
the error is reduced by a factor of $\kappa$ at each step until the sum 
of the significant coefficients in  {\rm (\ref{e41})}
is smaller than the final accuracy, 
which by {\rm (\ref{e42})}    
means that the same  property holds for the current 
Galerkin approximation.
\item[iv)] Although the sum in the right-hand 
side of {\rm (\ref{e40})} formally still contains 
unfinitely many coefficients,
it can be checked that this sum in fact  runs over a finite set, 
so that the adaptive strategy is  implementable.
\end{itemize}  } 
\end{rem}

Let us now compare this concept of adaptivity with the 
notion of adaptive information
explained above:

\begin{itemize}
\item 
{}From the discussion presented above, 
we  have seen that  adaptive wavelet schemes
are not  performed by  gaining more and more information from the right-hand
side $f$ in an adaptive fashion. Instead they use
the {\it residual} which depends on the right-hand side,
the operator, and the  domain.  Moreover, we see that the starting
index set $\Lambda$  is determined by the 
wavelet expansion of the right-hand side. That is,  $\Lambda$ is 
given by some kind of best $n$-term approximation of $f$,  which is assumed
to be available or to be  easily computable. In this sense, the adaptive
wavelet schemes require  {\it nonlinear information} about the problem.
 
\item 
In the wavelet setting, the benchmark for the performance 
is the approximation order of the best $n$-term approximation
of the solution, i.e., the numbers 
\begin{equation}   \label{ddd} 
\sup_{\| f\|_F \le 1} \, \sigma_n ( \A^{-1} f,\B)_H .  
\end{equation} 
It has been shown quite recently in 
\cite{CDD} that a judicious variant of the algorithm  
outlined above  gives rise to the
same order of approximation as best $n$-term approximation, 
while the  number of  arithmetic operations that are needed stays proportional
to the number of unknowns. 
Here the authors implicitly assume that certain subroutines for fast 
matrix-vector multiplications, approximations
of the right-hand sides and  for thresholding are available, and that all  
these routines have to realize a given approximation rate.
Moreover, it is assumed that the solution  $u$ is contained in some Besov 
space $B^\alpha_p(L_p(\Omega))$, and 
hence $F$ is a suitable subset of 
$\A (B^\alpha_p(L_p(\Omega)))$, i.e., 
the admissible class of right hand sides depends on the operator
$\A$. 
Observe that, for given $F$ and $\B$, the numbers 
$e_{n,C}^\non (S,F, H)$ might be much 
smaller than the numbers in \eqref{ddd} since it is, in general, 
not clear whether a wavelet basis is optimal. 

\item  
The performance of an adaptive scheme 
is not compared with  an {\it arbitrary linear} scheme.
The reason for that is simple,  and has already
been explained earlier. It is indeed true that linear 
approximation often produces the same order as 
nonlinear (best $n$-term) approximations, see 
Theorem \ref{t2} and Theorem \ref{t3}.
However, for nonregular problems,  it would be necessary to precompute
the optimal basis $S(g_i)$  in advance,  
which is mostly too expensive and should be avoided 
in practice, see \cite{DNS1} for further details. 
One usually compares
adaptive schemes with {\it uniform} methods for then a precomputation is 
not necessary. Therefore the use of an adaptive wavelet scheme is justified
if it performs better than any uniform scheme. It is known
that the order of approximation of uniform schemes is determined
by the Sobolev regularity $H^t(\Omega)$ of the object 
we want to approximate  whereas
the approximation order of best $n$-term approximation depends
on the regularity in the  specific Besov scale $B^t_{\tau}(L_{\tau}(\Omega))$,
where 
$$
\frac{1}{\tau}=\frac{t-s}{d}+\frac{1}{2},
$$
see {\rm \cite{DDD, DV98}} for details.  
Therefore adaptive schemes
are justified if the Besov regularity of the exact solution is
higher than its Sobolev regularity. 
For elliptic  boundary value problems, there
exist now many results in this direction, 
see, e.g., \cite{Da1,D1,D2,D4,DD}.   

\item 
In approximation theory, an approximation scheme that comes
from a sequence of linear spaces that are uniformly refined
is  also called {\it linear approximation scheme},  which sometimes
causes misunderstandings because these schemes are only special
cases of the linear schemes  considered, e.g., in Theorem~\ref{t3}. 
To avoid
this confusion, we used the term uniform methods instead of linear methods. 
\end{itemize}

\begin{rem}  
In this paper we study the complexity of solving elliptic 
partial differential equations. We only deal with the deterministic 
setting. The randomized setting, where also the use of random numbers
is allowed, is studied by Heinrich~{\rm \cite{H06}}. 
The complexity of solving elliptic PDE in the quantum model 
of computation (where one can use a certain nonclassical randomness) 
is studied in {\rm \cite{H06a}}. 
\end{rem} 


\section{Proofs}  \label{proofs}  



\subsection{Properties of  Widths}   
\label{sub4.1}


{\bf Proof of Lemma} \ref{bernstein}.
{\em Step 1.}
Part (i) is proved in \cite{TWW88} for the case where $F$ 
is normed. The general case is similar. 
\\
{\em Step 2.}
To prove part (ii),  we 
assume that $S(\{ \Vert f \Vert_F \le 1 \} )$ contains an 
$(n+1)$-dimensional ball $B \subset H$ of radius $r$ 
and  that $N_n : F \to \R^n$ is continuous. 
Since $S^{-1}(B)$ is an $(n+1)$-dimensional bounded symmetric 
neighborhood of 0, it follows from the Borsuk Antipodality Theorem, 
see \cite[paragraph 4]{D85},
that there exists an $f \in \partial S^{-1} (B)$ with 
$N_n (f)=N_n (-f)$ and hence 
$$
S_n(f) = \phi_n (N_n (f)) = \phi_n  (N_n (-f))= S_n(-f)
$$
for any mapping $\phi_n  : \R^n \to G$. 
Observe that $\Vert f \Vert_F=1$. 
Because of $\Vert S(f) - S(-f) \Vert =2r$ and 
$S_n(f)=S_n(-f)$ we obtain that the maximal error of $S_n$ 
on $\{ \pm f \}$ is at least~$r$. 
This proves 
\[
b_n (S, F , H) \le e_n^\cont (S, F , H)\, . 
\]
Since we  did not use the continuity of $\varphi_n$ also
$b_n (S, F , H) \le d^n (S, F , H)\, $ follows. \eproof
{~}\\

\noindent
{\bf Proof of Lemma} \ref{basic}. 
{\em Step 1.} Proof of (i). A corresponding assertion with $X$ 
and $Y$ normed linear spaces has been proved in \cite{DD96}.
This proof carries over without changes.
\\
{\em Step 2.} Proof of (\ref{e25}).
Let $\B= \{h_1,h_2,\ldots \, \}$ 
be a Riesz basis of $G$ with Riesz constants $A,B>0$. Let this basis 
$\B$ 
and a corresponding mapping
$S_n$ be optimal with respect to $I,F,G$ (up to some 
$\varepsilon>0$ if necessary). Then the image of $\B$ under the mapping
$S$ is a Riesz basis of $H$ with 
Riesz constants $A' = A/\|S^{-1}\|$ and $B'= B \, \| S \|$.
{}From
\[
\| \, S f - (S \circ S_n) f\, \|_H \le 
\| S \|\, \| \, f - S_n (f)\, \|_G 
\]
it follows that 
\[
e^\non_{n,C\, \| S^{-1}\| \, \| S \| } (S,F,H) \le \| S\| \,  
e^\non_{n,C} (I,F,G) \, . 
\]
Replacing $C$ by $C/(\| S^{-1}\| \, \| S \|)$, 
the right-hand side in
(\ref{e25}) follows. \\
Now, let $\B \subset H$ be a Riesz basis with Riesz constants $A,B>0$.
Let $\B$ and a corresponding $S_n$ be optimal with respect to 
$S,F,H$ (again up to  some $\varepsilon>0$ if necessary). From
\[
\| \, I f - (S^{-1} \circ S_n) f\, \|_G \le 
\| S^{-1} \|\, \| \, Sf - S_n (f)\, \|_H 
\]
it follows that 
\[
e^\non_{n,C\, \| S^{-1}\| \, \| S \| } (I,F,G) \le \| S^{-1}\| \,  
e^\non_{n,C} (S,F,H) \, . 
\]
The proof of (\ref{e24}) follows from (\ref{mult}).
\eproof
{~}\\

\noindent
Next we turn to the proof of Theorem \ref{t1}. 
It is convenient for us to start with a simplified situation.
For this we assume that $K \subset H$ is compact. 
We define
\begin{equation}   \label{ab01} 
e_{n,C}^\non (K, H) = \inf_{\B \in \B_C} \sup_{u \in K} \sigma(u, \B)
\end{equation} 
and
\begin{equation}     \label{ab02} 
e_n^{\rm cont} (K, H) = \inf_{N_n, \phi_n} \sup_{u \in K} \Vert \phi_n
(N_n(u)) - u \Vert ,
\end{equation} 
where the infimum runs over all continuous mappings 
$\phi_n : \R^n \to H$ and $N_n : K \to \R^n$. 
We prove the following result.

\begin{proposition}\label{vergleich*}     
Let $K \subset H$ be compact. Then
\begin{equation}   \label{ab03} 
e_{4n+1}^\cont (K,H) \le 2 C \, e_{n,C}^\non (K, H) .
\end{equation} 
\end{proposition} 

\begin{proof} 
Let $\B \in  \B_C$ be given. 
Since $K$ is compact, we only need finitely many elements 
of $\B$,  in the sense that 
\begin{equation} 
\sup_{u \in K} \Vert u - L_N (u) \Vert \le \e
\end{equation} 
for 
\begin{equation}   \label{ab05} 
L_N (u) = \sum_{j=1}^N a_j h_j .
\end{equation} 
Here $L_N$ is the orthogonal projection onto the space 
that is generated by $h_1, \dots , h_N$. The functionals $a_j$ 
are linear and continuous. Moreover, we know that 
\begin{equation}   \label{ab06} 
A \left( \sum_{j=1}^N |\a_j |^2 \right)^{1/2}   \le 
\Vert \sum_{j=1}^N \a_j h_j \Vert \le
B \left( \sum_{j=1}^N |\a_j |^2 \right)^{1/2}  
\end{equation} 
with $B/A \le C$. 
We may assume that $A=1$. 
For a suitable $\B \in \B_C$ we obtain 
\begin{equation}  \label{ab07} 
\sup_{u \in K} \biggl\| 
\sum_{k=1}^n c_k \,  h_{i_k}  - L_N (u) \biggr\| \le
e_{n,C}^\non (K,H) +\e . 
\end{equation} 
Let $\beta >0$. We define a modification 
of $L_N$ by
\begin{equation} 
L_N^* (u) = \sum_{j=1}^N a_j^*  h_j 
\end{equation} 
where
$a_j^* = a_j$ if $|a_j| \ge 2 \beta$ and
$a_j^* =0$ if $|a_j| \le \beta$. 
To make the $a_j^*$ continuous we define 
$$
a_j^* = 2 \,  \sgn (a_j) \cdot (|a_j| -\beta ) 
$$
for $|a_j| \in (\beta, 2\beta)$. 
We prove certain statements about $L_N^*$ and denote the best 
$n$-term approximation of $u$ by $u_n$. 

Assume that for $u \in K$, there are $m>n$ of the $a_j$, see
\eqref{ab05}, such that $|a_j| \ge \beta$. 
Then we obtain 
$$
\Vert u_n - L_N(u) \Vert \ge (m-n)^{1/2} \beta 
$$
and with \eqref{ab07} we obtain 
\begin{equation}  \label{ab08} 
m-n \le \frac{1}{\beta^2} ( e_{n,C}^\non (K, H) +\e)^2 .
\end{equation} 
Now we consider the sum $\sum_{|a_j|<\beta} a_j^2$ for $u \in K$. 
We distinguish between those $j$ that are used for $u_n$ 
(there are only $n$ of those $j$) and the other indices and 
obtain 
$$
\sum_{|a_j|^2<\beta} a_j^2 \le n \beta^2 + (e_{n,C}^\non(K,H) +\e)^2 .
$$
Now we are ready to estimate $\Vert L_N^*(u) - L_N(u)\Vert$ 
for $u \in K$. Observe that $|a_j^* - a_j| \le \beta$ 
for any $j$. We obtain 
$$
\Vert L^*_N (u) - L_N(u) \Vert \le 
B ( m\beta^2 + n \beta^2 + (e_{n,C}^\non(K,H) +\e)^2 )^{1/2} .
$$
Using the estimate \eqref{ab08} for $m$,  we obtain 
$$
\Vert L^*_N (u) - L_N(u) \Vert \le 
B ( 2  n \beta^2 + 2 (e_{n,C}^\non(K,H) + \e)^2)^{1/2} . 
$$
Now we define $\beta$ by
$$
n\beta^2 = ( e_{n,C}^\non (K,H) +\e)^2 
$$
and obtain the final error estimate (where we replace, for general $A$,
the number $B$ by $B/A$) 
$$
\Vert L_N^* (u) - L_N(u) \Vert \le 
\frac{2B}{A} \,  (e_{n,C}^\non(K,H) +\e) .
$$
In addition we obtain 
$$
m \le 2n 
$$
and therefore $L_N^*$ yields a
continuous $2n$-term approximation of 
$u \in K$ with error at most
$$
\sup_{u \in K} \Vert L_N^* (u) - u \Vert \le 
\frac{2B}{A} \,  (e_{n,C}^\non(K,H) +\e) + \e .
$$
The mapping $L_N^*$ is continuous and the image is 
a complex of dimension $2n$, see, e.g.,  \cite{DKLT93}. 
Hence we have an upper bound for the so-called \emph{Aleksandrov 
widths}, see \cite{DKLT93} and \cite{St74}.  
By the famous theorem of N\"obeling, 
any such mapping can be factorized as
$L_N^*=\phi_{4n+1} \circ N_{4n+1}$ where 
$N_{4n+1} : K \to \R^{4n+1} $ and 
$\phi_{4n+1} : \R^{4n+1} \to H$ are continuous. 
Hence the result is proved. 
\end{proof} 

\noindent
{\bf Proof of Theorem} \ref{t1}.
The unit ball of $F $  is a compact subset of $G$ by assumption. 
{}From Proposition \ref{vergleich*},  we derive that 
\[
e_{4n+1}^\cont (I,F,G) \le 2 C \, e_{n,C}^\non (I,F,G)\, .
\]
Next we apply Lemma \ref{basic}(ii),  and obtain
\[
e^\cont_n (S,F,H) \le  \| S \| \, 
e^\cont_n (I,F,G) , 
\]
as well as
\[
e^\non_{n,C} (I,F,G) \le \| S^{-1}\| \,  
e^\non_{n,C/(\| S^{-1}\| \, \| S \| )} (S,F,H) \, . 
\]
Combining these inequalities,   we are done. \eproof


\subsection{Widths of Embeddings of Weighted Sequence Spaces}  
\label{sequence}  


Having the wavelet characterization of Besov spaces in mind, cf. 
Subsections \ref{app3} and~\ref{app4}, we introduce the 
following scale of sequence spaces.

\begin{definition}
Let $0< p,q \le \infty$ and let $s \in \R$. 
Let $\nabla := (\nabla_j)_j$ be a sequence 
of subsets of finite cardinality of
the set $\{1,2,\ldots \,  , 2^d-1\} \times \Z^d$. 
We suppose that there exist $0< C_1 \le C_2$ 
and $J \in \Nb$ such that the 
cardinality $|\nabla_j|$ of $\nabla_j$ satisfies
\begin{equation}\label{301}
C_1 \le  2^{-jd} \, |\nabla_j| \le C_2 
\qquad \mbox{for all}\quad j \ge J \, . 
\end{equation} 
Then $b^s_{p,q}(\nabla)$,  where $ 0< q < \infty$, 
denotes the collection of all sequences
$a = (a_{j,\lambda})_{j,\lambda}$ of complex numbers such that
\begin{equation}\label{300}
\| \, a \, \|_{b^s_{p,q}}  := \left(  \sum_{j=0}^\infty  
2^{\displaystyle {j(s+d(\frac 12 -\frac 1p))q}} 
\bigg( \sum_{\lambda \in \nabla_j} 
| \, a_{j,\lambda} |^p \bigg)^{q/p}
\right)^{1/q} <\infty\,.
\end{equation}
For $q=\infty$, we  use the usual modification
\begin{equation} \label{300a}
\| \, a \, \|_{b^s_{p,\infty}}:= \sup_{j=1,2,\ldots}
2^{\displaystyle{j(s+d(\frac{1}{2}-\frac{1}{p}))}} 
\left(\sum_{\lambda \in \nabla_j}|a_{j,\lambda}|^p\right)^{1/p}
<\infty.
\end{equation}
If there is no danger of confusion we shall 
write $b^s_{p,q}$ instead of $b^s_{p,q}(\nabla)$. 
\end{definition}

\begin{rem}\label{iso}
In what follows,  we shall let $e_{j,\lambda}$ denote  the 
elements of the canonical orthonormal basis of $b_{2,2}^0$.
Let $\sigma \in \R$.
It is obvious that the linear mapping $L_\sigma$ defined by
\[
L_\sigma \, e_{j,\lambda} := 2^{-\sigma j}\, e_{j,\lambda} 
\qquad \mbox{for all}\quad j,\, \lambda\, , 
\]
extends to an isomorphism from $b^s_{p,q} $ 
onto $b^{s+\sigma}_{p,q}$ (simultaneously for all $s,p,q$)
with $\| \, L_\sigma \, \| =1$.
\end{rem}

\noindent
In the framework of these sequence spaces it is very easy to prove 
embedding theorems, cf. \cite{Le99}.

\begin{lemma}\label{einbettung}
Let $0 < p_0,p_1, q_0, q_1 \le \infty$, $s \in \R $, and $t\ge 0$.\\
{\rm (i)} The embedding  
\[
b^{s+t}_{p_0,q_0}(\nabla)\hookrightarrow b^{s}_{p_1,q_1}(\nabla)
\]
exists (as a set theoretic inclusion) if and only if it is 
continuous if and only if either 
\begin{equation}\label{einb*}
t >  d \, \biggl(\frac{1}{p_0} - \frac{1}{p_1}\biggr)_+
\end{equation} 
or 
\[
t =  d \, \biggl(\frac{1}{p_0} - \frac{1}{p_1} \biggr)_+ 
\qquad \mbox{and}\qquad q_0 \le q_1\, .
\]
{\rm (ii)} The embedding  
\[
b^{s+t}_{p_0,q_0}(\nabla)\hookrightarrow b^{s}_{p_1,q_1}(\nabla)
\]
is compact if and only if {\rm (\ref{einb*})} holds. 
\end{lemma}

\noindent
The main result of this subsection consists in the following:

\begin{thm}   \label{th3} 
Let $0 < p, p_0,p_1 \le \infty$, $0 < q, q_0,q_1 \le \infty$, 
and $s \in \R $.\\
{\rm (i)} Suppose that
\begin{equation}\label{00}
t> d\, \Big(\frac{1}{p} - \frac 12 \Big)_+
\end{equation}
holds. Then, for any $C\ge 1$,  we have 
\[
e_{n,C}^\non (I, b^{s+t}_{p,q},   b^s_{2,2} )   \asymp 
n^{t/d} \, .  
\]
{\rm (ii)} Suppose that {\rm (\ref{00})} holds.
Then we have 
\[
e_{n}^\lin (I, b^{s+t}_{p,q},   b^s_{2,2} )   \asymp \left\{
\begin{array}{lll}
n^{-t/d} & \qquad & \mbox{if}\quad 2 \le p \le \infty  , \\
n^{-t/d+1/p-1/2}  & 
\qquad & \mbox{if}\quad 0< p <2 . 
\end{array}\right.
\]
{\rm (iii)} Suppose that {\rm (\ref{einb*})} holds. Then we have 
\[
e_n^\cont (I, b^{s+t}_{p_0,q_0},   b^s_{p_1,q_1} )   \asymp 
n^{-t/d} \, .
\]
\end{thm} 

\begin{rem}
\label{Limesfall} 
In part (i) there is an interesting limiting case. Suppose
$0 < p < 2$ and $t= d (1/p-1/2)$.
Then the embedding $b^{s+t}_{p,p} \hookrightarrow b^s_{2,2}$ exists, 
cf. Lemma \ref{einbettung}, and 
\[
\left(\sum_{n=1}^\infty \Big[ n^{t/d} \, \sigma_n (a,\B)_{b^s_{2,2}}\Big]^p 
\, \frac 1n  \right)^{1/p} < \infty 
\quad \mbox{if and only if} \quad a \in  b^{s+t}_{p,p} \, .
\]
In view of Lemma \ref{einbettung}(ii),  
this shows that $\lim_{n\to \infty} e^\non_{n,C} (S,F,H) =0$ 
does not imply compactness of $S$.
\end{rem}

\noindent
The proof of Theorem \ref{th3} requires some preparations.
It will be given in Subsections \ref{aaa}--\ref{aac}.


\subsubsection{The Bernstein  Widths of the Identity Operator}   


We concentrate on the estimate from below.
For later use we treat  a more general situation.

\begin{lemma}\label{bernst}
Let $0 < p_0,p_1, q_0, q_1 \le \infty$, $s \in \R $ and $t>0$
such that {\rm (\ref{einb*})} holds. 
Then there exists a positive constant $c$
such that
\begin{equation}
b_n   (I, b^{s+t}_{p_0,q_0},  
b^s_{p_1,q_1}) \ge c\, \left\{
\begin{array}{lll}
n^{-t/d} & \qquad & \mbox{if} 
\quad 0 < p_0 \le p_1 \le \infty \, , \\
n^{- t/d+1/p_0 -1/p_1}  && \mbox{if} 
\quad 0 <  p_1 < p_0 \le \infty\, .
\end{array}
\right.
\end{equation}
holds for all $n$.
\end{lemma}

\begin{proof}
The Bernstein numbers are monotonic in $n$. 
So it will be enough to prove the assertion
for sufficiently large $n$.    
Consequently, we may assume that there is a natural 
number $N\ge J$, as well as positive constants $c_1$ and $c_2$,  such that
\[
c_1 \, 2^{Nd} \le n \le c_2 \,2^{Nd} \, .
\]
{\em Step 1.}
Let $0 < p_0 \le p_1$.
Using H\"older's inequality we find
\begin{eqnarray*}
\| \, \sum_{\lambda \in \nabla_N} \, b_\lambda \, e_{N,\lambda}  
|b^{s+t}_{p_0,q_0} \|
& = &  
 2^{N(s + t + d/2-d/p_0 )} \, 
\Big(\sum_{\lambda \in \nabla_N} |b_\lambda |^{p_0}\Big)^{1/p_0}\\
& \le &   2^{N(s + t  + d/2      -d/p_0)} \, 
|\nabla_N|^{1/p_0 - 1/p_1}
\Big(\sum_{\lambda \in \nabla_N} |b_\lambda |^{p_1}\Big)^{1/p_1}\\
& \le & C_2\,  2^{Nt} \,  \| \, \sum_{\lambda \in \nabla_N} \, 
b_\lambda \, e_{N,\lambda} \, |b^{s}_{p_1,q_1}\|
\\
& \le & c_3\,  n^{t/d}\, 
\| \, \sum_{\lambda \in M_N} \, b_\lambda \, e_{N,\lambda}\,   
|b^{s}_{p_1,q_1}\| \, ,
\end{eqnarray*}
where $C_2$ corresponds to (\ref{301}). 
Consequently, the  unit ball in $b^{s+t}_{p_0,q_0}$ contains the
$n$-dimensional ball 
(spanned by the   vectors $e_{N,\lambda}$, 
$\lambda \in \nabla_{N}$) with  radius $c_3^{-1} \, n^{-t/d}$.
This proves
\[
b_n   (I, b^{s+t}_{p_0,q_0}, b^s_{p_1,q_1} ) \ge c\,  n^{-t/d}
\]
for some positive constant $c$ independent of $n$.
\\
{\em Step 2.}
If $p_0 > p_1$, then H\"older's inequality 
(used in the second line of the estimate in Step 1) 
will be replaced by the monotonicity of the $\ell_r$-norms and  we obtain
\begin{eqnarray*}
\| \, \sum_{\lambda \in \nabla_N} \, b_\lambda \, e_{N,\lambda} \, 
|b^{s+t}_{p_0,q_0}\|
& = &  \,  2^{N(s + t + d/2 -d/p_0)} \, 
\Big(\sum_{\lambda \in \nabla_N} |b_\lambda |^{p_0}\Big)^{1/p_0}
\\
& \le &   2^{N(s + t + d/2  -d/p_0)} \, 
\Big(\sum_{\lambda \in \nabla_N} |b_\lambda |^{p_1}\Big)^{1/p_1}
\\
& \le & c_5\,  2^{N( t +d/p_1 - d/p_0})\, 
\biggl\| \, \sum_{\lambda \in \nabla_N} \, 
b_\lambda \,  e_{N,\lambda}\, |b^{s}_{p_1,q_1} \biggr\| \, .
\end{eqnarray*}
This time  the  unit ball in $b^{s+t}_{p_0,q_0}$ contains the
$n$-dimensional ball  with radius 
$$
c_5^{-1}\,  2^{-N( t + d/p_1 - d/p_0)}. 
$$
This proves our claims.
\end{proof}

\begin{rem}\label{referenz}
In the one-dimensional periodic situation, estimates 
of the Bernstein numbers from above are also known, 
due to Tsarkov and Maiorov, cf. {\rm \cite[Thm.~12, p.~194]{Ti90}.}
Let $1\le p \le \infty$ and $s>0$. By
$\mathring{W}^s_p $ we denote the collection 
of all $2\pi$-periodic functions $f$
with Weyl derivative of order $s$ belonging to $L_p (\T)$ and satisfying
$\int_{-\pi} ^\pi f (x)\, dx =0$. Then 
\[
b_n   (I, \mathring{W}^{t}_{p_0},  L_{p_1}) \asymp \, \left\{
\begin{array}{lll}
n^{-t} & \qquad & \mbox{if} 
\quad 1\le  p_0 \le p_1 \le \infty \quad \mbox{or} \\
&& \hspace{1cm} 1 \le p_1 \le p_0 \le 2 \quad \mbox{and}\quad t>0\, , \\
n^{- t  + 1/p_0 -1/p_1} && \mbox{if} 
\quad 2 \le   p_1 < p_0 \le \infty \quad \mbox{and}
\quad t>1/p_0\, , \\
n^{- t  + 1/p_0 -1/2} && \mbox{if} 
\quad 1 \le p_1 \le 2 \le  p_0 \le \infty \quad 
\mbox{and}\quad t>1/p_0\, . \\
\end{array}
\right.
\] 
This should be compared with Lemma \ref{bernst} for $s=0$ and  $d=1$.
\end{rem}


\subsubsection{Best $m$-Term Approximation in the 
Framework of Sequence Spaces}   
\label{aaa}

We prepare the proof of part (i) of Theorem \ref{th3}.
Also,  here we treat a more general situation. 
Let $\B$ denote the canonical basis
$(e_{j,\lambda})_{j,\lambda}$ in $b^0_{2,2} (\nabla)$.
Then our aim in this subsection consists in  
a characterization of the behaviour of the best 
$m$-term approximation 
of a given element $a \in b^{s+t}_{p_0,q_0}$ with respect to $\B$.\\
The main result of this subsection reads as follows:

\begin{thm}   \label{t5*} 
Let $0 < p_0, p_1, q_0, q_1 \le \infty $, $s\in \R $ and $t>0$
such that {\rm (\ref{einb*})} holds.
Then we have 
\begin{equation}\label{claim2*}
\sup \Big\{  \sigma_n (a,\B)_{b^{s}_{p_1,q_1}} \, :  
\quad \|\, a\|_{b^{s+t}_{p_0,q_0}} \le 1\Big\} 
 \asymp  n^{-t/d}\,  . 
\end{equation}
\end{thm} 

\noindent
We start with some preparations.
Let $U $ denote the unit ball in $b^{s+t}_{p_0,\infty} $. Then 
\[ 
a = \sum_{j=0}^\infty \sum_{\lambda \in \nabla_j} 
a_{j,\lambda} \, {e}_{j,\lambda}  
\qquad \mbox{and}\qquad 
\sup_{j=0,1, \ldots \, } \, 2^{j(s + t +  d(1/2-1/p_0))} 
\Big(\sum_{\lambda \in \nabla_j} | a_{j,\lambda}|^{p_0}\Big)^{1/p_0} 
\le 1 \,  .
\] 
The following lemma will be of some use:

\begin{lemma} \label{mann*} 
Let $0 < p_0 \le p_1$ and suppose that 
\begin{equation}\label{einb2*}
t > d \, \Big( \frac{1}{p_0} - \frac{1}{p_1}\Big)\, .
\end{equation}
For all  $a \in U$ and all $n\ge 1$ there exists 
a natural number $ K := K(a,n)$  such that
\[
  \Big\| a - \sum^K_{j=0} 
\sum_{\lambda \in \nabla_j} a_{j,\lambda} \, e_{j,\lambda} 
\, \Big|  b^s_{p_1,q_1} \Big\| \ \le \ 
n^{-t/d}
\]
holds.
\end{lemma}

\begin{proof} We define
$$  T_j\,:= \ 
\sum_{\lambda \in \nabla_j} a_{j,\lambda} \, e_{j,\lambda}\, , 
\qquad j = 0 , 1 \ldots \, .
$$
Then one has
\[
 a - 
\sum^K_{j=0} \sum_{\lambda \in \nabla_j} a_{j,\lambda} \, e_{j,\lambda} \ 
= \ \sum_{j>K} T_j.
\]
Since of $0< p_0  \le p_1 \le \infty $,  the monotonicity 
of the $\ell_q$-norms  and $a \in U$ lead to 
\begin{eqnarray*}
  \|\, T_j\, |  b^s_{p_1,q_1} \| & \le &  
\,   2^{j  (s + d/2 - d/p_1 )} 
\Big(\sum_{\lambda \in \nabla_j} |a_{j,\lambda}|^{p_0}  \Big)^{1/p_0}\\
    & \le &  \, 2^{-j(t + d(1/p_0 - 1/p_1))}. 
\end{eqnarray*}
Let $u=\min (1,p_1,q_1)$.
Consequently,  using (\ref{einb2*}) and chosing $K$  large enough,   we find
\begin{eqnarray*}
  \Big\| \sum_{j\ge K} T_j \, \Big|  b^s_{p_1,q_1} \Big\|^u & 
  \le & \sum_{j\ge K} \|\, T_j \, |  b^s_{p_1,q_1} \|^u \ \le \ 
\sum_{j\ge K} 2^{-ju{\big[t + d(1/p_0 
- 1/p_1)\big]}}  \\
   & \le & C_1 \, 2^{-Ku(t + d(1/p_0-1/p_1))}
   \leq n^{-tu/d} . 
\end{eqnarray*}
This proves the claim. 
\end{proof}

\noindent
The basic step in deriving an  upper estimate of 
$\sigma_n (a,\B)$   
is the following proposition. 
Again $U $ denotes the unit ball in $b^{s+t}_{p_0,\infty} $. 

\begin{proposition} \label{proposi*}  
Let $0 < p_0 \le p_1 \le \infty$.
Let  $ a \in U$, $n \in \Nb$,  and let  $K=K(a,n)$ 
be as in Lemma \ref{mann*}. 
Then there exists an  approximation
\begin{equation}\label{Sdef1}
S_n a \ := \ \sum_{j=0}^K \sum_{\lambda \in \nabla_j} a^*_{j,\lambda} \, 
e_{j,\lambda}
\end{equation}
of $a$, which satisfies the following: 
\begin{itemize}
\item[i)] The coefficients $a_{j,\lambda}^*$ depend continuously on $a$.
\item[ii)] The number of nonvanishing entries is bounded by $c\cdot n$.
\item[iii)] \ $\| \, a - S_n a\, |  b^s_{p_1,q_1}\| \ 
\le \ c\,n^{-t/d}\,,\ \ n = 1,2,\ldots\,.$
\end{itemize}
Here $c$ can be chosen independent of $a$ and $n$.
\end{proposition}

\begin{proof}
Observe that it will be enough to prove the claim 
for natural numbers $n=2^{Nd}$, where  $N\in \Nb$.
We define
\begin{eqnarray}
\delta  & := & \frac{t - d\, \big(1/p_0 - 1/p_1\big)}
{2\big(1/p_0 - 1/p_1\big)} \, ,
\nonumber
\\
\varepsilon_j & := & \left\{ \begin{array}{lll}
    0 & \qquad & \mbox{if}\quad  1\le j\le N \\[4mm]
n^{-1/p_0} 2^{-jd (1/2-1/p_0)} 2^{-jt} 
2^{(j-N)\delta/p_0}
&& \mbox{if}\quad j > N \, , 
\end{array} \right. 
\label{Lambda*}
\\  \nonumber 
~ && \\ 
\label{Lambda}
\Lambda^*_{j}  & := &  \Big\{\lambda \in \nabla_j\,: \: 
\, |a_{j,\lambda}| \, 2^{s j} \ge \varepsilon_j\Big\}
\, , \qquad j=0,1,\ldots \, \, . 
\end{eqnarray}
Then, if $j >N$, 
\begin{eqnarray} \label{anzahl*}
|\Lambda^*_{j}| & = & \sum_{\lambda \in \Lambda^*_j} 1 \ 
\le \ \sum_{\lambda \in \Lambda^*_{j}} 
2^{j s p_0} \, \frac{|a_{j,\lambda}|^{p_0}}{\varepsilon^{p_0}_j} \\
& \le & \sum_{\lambda \in \nabla_j}  n \, 
2^{jd (1/2 - 1/p_0)p_0}  
2^{jtp_0} 2^{-(j-N)\delta} 2^{jsp_0} |a_{j,\lambda}|^{p_0} 
\nonumber
\\
& = & n\,2^{-(j-N)\delta} \, 
\sum_{\lambda \in \nabla_j} \, 
 2^{j(s + t + d(1/2 -1/p_0))p_0} \, |a_{j,\lambda}|^{p_0}
\nonumber
\\ 
& \le &  \ n\,2^{-(j-N)\delta}\, 
\| \, a \, |b^{s+t}_{p_0,\infty} \|^{p_0} 
\nonumber
\\
& \le &    n\,2^{-(j-N)\delta}\, .
\nonumber
\end{eqnarray}
Now a typical method to approximate $a$  would be to choose  
\ $a^*_{j,\lambda} = a_{j,\lambda}\,, \ j\in \Lambda^*_j$ 
and zero otherwise. 
However, this  selection does not depend continuously on $a$. 
Therefore we use the following variant.
Let  $g_j$ denote the following  piecewise  linear and odd function,
\begin{equation}\label{lambda2*}
  g_j(x) \ := \ \left\{ \begin{array}{lll} 
   0 &\qquad & \mbox{if}\quad  0\le x\le 2^{-js}\varepsilon_j \, , 
\\[2mm]
 x && \mbox{if}\quad x\ge 2\cdot 2^{-js}\varepsilon_j \, ,
\\[2mm]
   \mbox{linear} && \mbox{if} \quad x\in (2^{-js} \varepsilon_j, 
2\cdot 2^{-js}\varepsilon_j)\, . \end{array} \right.  
\end{equation}
Then we set 
\begin{equation}\label{lambda3*}
  a^*_{j,\lambda} \ := \ g_j (a_{j,\lambda})
\end{equation}
and consider the associated approximation (\ref{Sdef1}).
Let us prove that $S_n$ will do the job. 
\\
{\em Step 1.} We shall prove (i). 
Observe
\[
\Big|\, \bigcup_{j=0}^K \Lambda_j^* \, \Big| \le c_1 \, 
\sum_{j=0}^N \,  2^{jd} 
+  \, \sum_{j=N+1}^K \,  n \, 2^{-(j-N)\delta} \le c_2  \, n\, , 
\]
cf. (\ref{anzahl*}). 
The constant $c_2$ is independent of $a,K$, and $n$.
This proves (i) and (ii).\\
{\em Step 2.} Proof of (iii). We have
\[
a - S_n a \ = \ a - \sum^K_{j=0} 
\sum_{\lambda \in \nabla_j} a_{j,\lambda}\,  e_{j,\lambda}\, +\, 
\sum^K_{j=0} T^*_j  \ =: \ \Sigma_1 + \Sigma_2, 
\]
where
\[
T^*_j \ = \ \sum_{\lambda \in \nabla_j} 
\big( a_{j,\lambda} - a^*_{j,\lambda} \big) \, e_{j,\lambda} . 
\]
{}From Lemma \ref{mann*},  
we can conclude that \ $\|\, \Sigma_1 \, |  b^s_{p_1,q_1}\| 
\le n^{-t/d}$ 
for $K$ large enough. 
Therefore it remains to estimate $\|\, T^*_j \, |  b^s_{p_1,q_1}\|$.
Since  
 $|g_j(x) - x| \ \le \ |x|$ and 
$a^*_{j,\lambda} \ = \ a_{j,\lambda}\quad  \mbox{for} 
\quad |a_{j,\lambda}| \ \ge \ 2 \varepsilon_j 2^{-js}$, we obtain
\begin{eqnarray*}
|\, a_{j,\lambda} - a^*_{j,\lambda}\, |^{p_1} \ &
\le & \ |a_{j,\lambda}|^{p_1} \\
   & \le &  \ |a_{j,\lambda}|^{p_0} |a_{j,\lambda}|^{p_1-p_0} \\
&\le & \ |a_{j,\lambda}|^{p_0} 
(2\varepsilon_j)^{p_1-p_0} 2^{-js(p_1-p_0)} \, .
 \end{eqnarray*}
This will be used to estimate the norm of $T_j^*$ as follows:
\begin{eqnarray*}
\|\, T^*_j\, |  b^s_{p_1,q_1} \| & = &  \,  2^{j(s + d(1/2 - 
1/p_1))}
\bigg( \sum_{k\in \nabla_j} |a_{j,\lambda} - a^*_{j,\lambda}|^{p_1} 
\bigg)^{1/p_1} \\
 & \le & c_1 \, 2^{j d(1/2 - 1/p_1)} \, 2^{jsp_0/p_1}
\varepsilon_j^{1-p_0/p_1} 
\bigg( \sum_{k \in \nabla_j}  |a_{j,\lambda}|^{p_0}  \bigg)^{1/p_1} \\
     & \le & c_1 \, 
\varepsilon^{1-p_0/p_1}_j 
2^{j d/2} 2^{-jtp_0/p_1} 2^{-jd p_0/(2p_1)} 
\bigg( \sum_{\lambda \in \nabla_j} 2^{j(s + t + d(1/2  
- 1/p_0))p_0} 
|a_{j,\lambda}|^{p_0}  \bigg)^{1/p_1} \\
& \le & c_2 \, 
\varepsilon^{1-p_0/p_1}_j 
2^{-j (t + d/2 -dp_1/(2 p_0))p_0/p_1}\, 
\, \| \, a\, |b^{s+t}_{p_0,\infty} \|^{p_0/p_1} \\
 & \le & c_2 \, \varepsilon^{1-p_0/p_1}_j 
2^{-j (t + d/2 -dp_1/(2 p_0))p_0/p_1}\, , 
\end{eqnarray*}
where again $c_2$ does not depend on $a$ and $n$.
For $j>N$ we continue by employing the concrete value of $\varepsilon_j$
and obtain
\begin{eqnarray*}
\|\, T^*_j\, |  b^s_{p_1,q_1}\| & \le & c_2 \, 
\Big(n^{-1/p_0} 2^{-jd (1/2-1/p_0)} 2^{-jt} 
2^{(j-N)\delta/p_0}\Big)^{1-p_0/p_1} 
2^{-j (t + d/2 -dp_1/(2 p_0))p_0/p_1}
\\
 & = & c_2\,  n^{1/p_1 - 1/p_0} \,  
     2^{-N \delta( 1/p_0 - 1/p_1)} 
2^{-j(t - d(1/p_0 - 1/p_1) 
- \delta/p_0 + \delta/p_1) }\, .
\end{eqnarray*}
By construction $T_j^* =0$ if $j \le N$,  by definition, we have 
\[ 
t - d \Big(\frac{1}{p_0}- \frac{1}{p_1}\Big) 
> \delta \, \Big( \frac{1}{p_0} - \frac{1}{p_1}\Big)
\, . 
\]
Hence, with $u= \min (1,p_1,q_1)$, we have 
\begin{eqnarray*}
  \|\, \Sigma_2 \, |  b^s_{p_1, q_1}\|^u  & \le & c_2^u \,  
  \Big(n^{1/p_1 - 1/p_0} \,  
     2^{-N \delta (1/p_0 - 1/p_1)} \Big)^u
\sum^K_{j=N+1} 2^{-ju(t - d(1/p_0 - 1/p_1) 
- \delta/p_0 + \delta/p_1)}
\\
& \le & c_3 \, 
\Big(n^{1/p_1 - 1/p_0} \,  
     2^{-N \delta (1/p_0 - 1/p_1)} \Big)^u
\,  2^{-Nu(t - d(1/p_0 - 1/p_1) 
- \delta/p_0 + \delta/p_1)}
 \\
& = & c_3 \, 
\Big(n^{1/p_1 - 1/p_0} \Big)^u
\,  2^{-Nu(t - d(1/p_0 - 1/p_1))}
\, ,
\end{eqnarray*}
with $c_3$ independent of $K,n$ and $a$.
Recalling that   $2^{Nd}= n$, we end up with 
\[
  \|\, \Sigma_2\, |  b^s_{p_1,q_1}\| \ \le \ 
c_3  \, n^{-t/d} \,  .
\] 
This finishes the proof of Proposition \ref{proposi*}. 
\end{proof}

\noindent
For completeness and better reference we formulate 
the counterpart of Proposition~\ref{proposi*} in the case $p_0 \ge p_1$.

\begin{proposition} \label{proposi**}  
Let $0 < p_1 \le p_0 \le \infty$.
Let $ a \in U$
(the unit ball in $b^{s+t}_{p_0,\infty} $) 
and $2^{Nd} \le n \le 2^{(N+1)}d$. Then  the   approximation
\begin{equation}\label{Sdef2}
S_n a \ := \ \sum_{j=0}^N \sum_{\lambda \in \nabla_j} a_{j,\lambda} \, 
e_{j,\lambda}
\end{equation}
of $a$  satisfies the following: 
\begin{itemize}
\item[i)] The coefficients $a_{j,\lambda}$ depend continuously on $a$.
\item[ii)] The number of nonvanishing entries is bounded by $c\cdot n$.
\item[iii)] \ $\| \, a - S_n a\, |  b^s_{p_1,q_1}\| \ 
\le \ c\,n^{-t/d}\,,\ \ n = 1,2,\ldots\,.$
\end{itemize}
Here, $c$ can be chosen independent of $a$ and $n$.
\end{proposition}

\begin{proof}
The proof is elementary.
\end{proof}

\noindent
{\bf Proof of Theorem \ref{t5*}}.
The estimate from above follows from Propositions  
\ref{proposi*} and \ref{proposi**}, as well as 
the continuous embedding
$b^{s+t}_{p_0,q_0}  \hookrightarrow 
b^{s+t}_{p_0,\infty} $. 
For the estimate from below, 
it will be enough to consider $n=2^{Nd}$, where 
$N\ge J$ and $N \in \Nb$. 
Let $K$ be the smallest natural number such 
that $C_1 \, 2^{Kd} \ge 2$ (here $C_1$ is the same 
constant as in (\ref{301})).
Then
\[ 
n \le \frac{C_1 \, 2^{(N+K)d}}{2} \le  \frac 12 \,  |\nabla_{N+K}| \, .
\] 
Let $\Gamma \subset \nabla_{N+K}$ with $|\Gamma| = n$.
We define
\[ 
a = |\nabla_{N+K}|^{-1/p_0} \, 2^{-(N+K)(s+t + 
d(1/2 - 1/p_0))} \, \sum_{\lambda \in 
\nabla_{N+K}} e_{N+K,\lambda}\, .
\]
Consequently
$\| \, a \, \|_{b^{s+t}_{p_0,q_0}} = 1$ for any $q_0$.
Furthermore, we find
\begin{eqnarray*}
\| a - S_n a \|_{b^s_{p_1,q_1}} & \ge & \,  
\, \Big\| \sum_{\lambda \in \nabla_{N+K} \setminus \Gamma}
|\nabla_{N+K}|^{-1/p_0} \, 2^{-(N+K)(s+t + d(1/2 - 1/p_0))}\,  
e_{N+K,\lambda}\Big\|_{b^s_{p_1,q_1}} \\
& = &   |\nabla_{N+K}|^{-1/p_0} \, 
2^{-(N+K)(t + d(1/p_1 - 1/p_0))} 
|\nabla_{N+K} \setminus \Gamma|^{1/p_1}\\
& \ge & \frac{C_1^{1/p_1}}{2^{1/p_1}\, C_2^{1/p_0}} \, 2^{-(N+K)t}
\\
& = & \frac{C_1^{1/p_1}}{2^{1/p_1}\, C_2^{1/p_0}} \, 2^{-Kt}\, n^{-t/d}
\, ,
\end{eqnarray*}
(also $C_2$ has the same meaning as in (\ref{301})).
It is clear that an optimal $\Gamma$ with $|\Gamma|=n$
has to be a subset of $\nabla_{N+K}$.
This completes the proof of the estimate from below.
\eproof
{~}\\

\noindent
{\bf Proof of Theorem \ref{th3}(i)}.
The estimate from above is covered by Theorem \ref{t5*}; 
the estimate from below follows from Theorem \ref{t1} and 
Theorem \ref{th3}(iii).
\eproof


\begin{rem} 
Stepanets {\rm \cite{Ste}}  has investigated the quantities
$$ 
\sigma_n(a,B)_{b^s_{p_1,q_1}}
$$
for the specific case
$$ 
s=d \, \biggl(\frac{1}{p_1}-\frac{1}{2} \biggr), 
\quad \hbox{with} \quad  p_1=q_1.
$$
In this special case, the associated nonlinear withs related to 
quite general smoothness spaces are studied.
He proved explicit formulas from which the asymptotic 
behavior could be derived. 
\end{rem} 


\subsubsection{The Manifold Widths of the Identity}   
\label{aab}

{\bf Proof of Theorem \ref{th3}(iii)}.
Without loss of generality we may choose $s=0$, cf. 
Lemma \ref{basic}(ii) and Remark \ref{iso}.\\
{\em Step 1.} The estimate from above. In the case $p_1=q_1=2$ we may 
use Propositions \ref{vergleich*}, \ref{proposi*} 
and \ref{proposi**} to get the desired inequality.
However, for the general case we have to modify the argument.
We follow the arguments used in \cite{DKLT93}.
Let $U$ denote the unit ball in $b^{t}_{p_0,q_0} $.
As explained there Propositions  \ref{proposi*} and \ref{proposi**}
guarantee that 
\[
a^n (U,b^0_{p_1,q_1}) \le c \,  n^{-t/d}\, , 
\]
where $a^n $ denotes the Alexandroff-co-width, 
cf. \cite{DKLT93} for details.
But
\[
e^\cont_{2n+1} (U,b^0_{p_1,q_1}) \le a^n (U,b^0_{p_1,q_1})\, ,
\]
cf. \cite{DKLT93} and \cite{DD96}. 
Let us mention that in the  literature quoted the target 
space was always a normed linear space. 
But the arguments carry over to quasi-normed linear spaces.\\
{\em Step 2.} The estimate from below.
Lemmas \ref{bernstein} and \ref{bernst} yield the lower estimate in case
$0 < p_0 \le p_1 \le \infty$.\\
Now, let  $p_1 < p_0 \le \infty$. Let $\varepsilon >0$.
We consider the diagram
\begin{eqnarray*}
b^0_{p_1,q_1}  & \stackrel{\hbox{$I_3$}}{\longrightarrow} & 
b^{-d(1/p_1-1/p_0) - \varepsilon}_{p_0,\infty}  \\
I_2 \,  & \nwarrow \qquad  \nearrow & I_1 \\
& b^{t}_{p_0,q_0},  &
\end{eqnarray*}
where $I_1,I_2$ and $I_3$ are identity operators.
Then (\ref{Khod}) yields
\[
e_{2n}^{\rm cont} (I_1,b^{t}_{p_0,q_0} ,
b^{-d(1/p_1-1/p_0) - \varepsilon}_{p_0,\infty} ) 
 \le  
e_{n}^{\rm cont} (I_2,b^{t}_{p_0,q_0}, b^0_{p_1,q_1} )\, \,   
e_{n}^{\rm cont} (I_3,b^0_{p_1,q_1}, b^{-d(1/p_1
-1/p_0) - \varepsilon}_{p_0,\infty} )
\]
which implies that 
\[
c_1 \, n^{-t/d - 1/p_1+ 1/p_0 - \varepsilon/d} 
 \le c_2 \,  
e_{n}^{\rm cont} (I_2,b^{t}_{p_0,q_0} , b^0_{p_1,q_1})\, \,   
n^{- 1/p_1+1/p_0 - \varepsilon/d}
\]
for some positive $c_1$ and $c_2$ (independent of $n$), see
Lemmata \ref{bernst}, \ref{bernstein}, and Step 1. \eproof
{~}\\

\begin{rem} 
It is clear from the proof given above that the knowledge
of the Bernstein widths is not enough
to 
establish the estimate from below of 
$e^\cont_n$.
Here the multiplicativity of the numbers $e^\cont_n$, cf. 
{ \rm (\ref{Khod})}, is crucial.
This seems to be overlooked in {\rm \cite{DKLT93}}.
\end{rem}


\subsubsection{The Approximation Numbers  of the Identity}   
\label{aac}


{\bf Proof of Theorem \ref{th3}(ii)}. 
{\em Step 1.} Let $2\le p \le \infty$. From
Proposition \ref{proposi**} we obtain the estimate from above 
with $S_n$ given by (\ref{Sdef2}).
The estimate from below is covered by (\ref{claim2*}).\\
{\em Step 2.} Let $0< p < 2$. Without loss of generality we assume $s=0$.
Let $S_n$ be defined by  (\ref{Sdef2}).
The estimate from above is easily derived 
by using the monotonicity of the $\ell_r$-norms 
and $t+d(1/2-1/p)>0$:
\begin{eqnarray*}
\| \, a-S_n a\, |b_{2,2}^0\|^2 & \le &
\sum_{j=N+1}^\infty  \Big(\sum_{\lambda \in \nabla_j}\, 
|a_{j,\lambda}|^{p}\Big)^{2/p}\\
&\le & \Big(\sum_{j = N+1}^\infty  2^{-2j (t+d(1/2 -1/p))}\Big)
\bigg( \sup_{j\ge N+1} 2^{j(t+d(1/2 -1/p))} 
\Big( \sum_{\lambda \in \nabla_j}\, |a_{j,\lambda}|^{p}\Big)^{1/p}\bigg)^2\\
&\le & c\, 2^{-2N(t+d(1/2 -1/p))} \| \, a \, |b^t_{p,\infty}\|^2
\\
&\le & c\, \Big(n^{-t/d - 1/2 +1/p} \| \, 
a \, |b^t_{p,q}\|\Big)^2\, , 
\end{eqnarray*}
where $c$ does not depend on $n$ and $a$.
For the estimate from below,  we use the obvious 
fact that the optimal approximation 
of an element in a Hilbert space is given by the 
partial sum with respect to an orthonormal basis.
Hence, if $\widetilde{S}_n$ is a linear operator of rank 
at most $n$ then
\[
\| \, a - \widetilde{S}_n a|b_{0,0}\|
\ge \| \, a - {S}_n a|b_{0,0}\|\, ,
\]
where $S_n$ is defined by  (\ref{Sdef2}).
We put 
\[
a:= \sum_{j=0}^{N+1}\,  e_{j,\lambda_j}\, , 
\]
where $\lambda_j \in \nabla_j$ can be chosen arbitrarily.
Then
\[
\| \, a \, |b^t_{p,q}\| = \Big(\sum_{j=0}^{N+1} 2^{j(t+d(1/2 -1/p))q}
\Big)^{1/q} \ge  2^{N(t+d(1/2 - 1/p))}
\]
for some positive $c$ independent of $n$ and
\[
\| \, a - S_n a |b_{2,2}^0\| = 1\, .
\]
This implies
\[
\| \, I - S_n \, | b^t_{p,q}\| \ge \frac{1}{2^{N(t+d(1/2 -1/p))}}\, ,
\]
which finishes the proof of the lower bound.
\eproof

\begin{rem}\label{opt}
Notice that in any case, 
an order-optimal approximation is given by an appropriate partial sum,
see {\rm (\ref{Sdef2})}.
\end{rem}


\subsubsection{The Gelfand  Widths of the Identity}   \label{Gelfand}


What we will do here relies on a result of Gluskin \cite{Gl1,Gl2}
about the Gelfand widths 
of the embedding $\ell_p^m \to \ell_2^m$ which we now recall.  
Let $1/p + 1/p'=1$.
For all natural numbers $m$ and $n$, where  $n\le m$, 
it holds that 
\begin{equation}\label{Gluskin}
d^n (I, \ell_p^m, \ell_2^m) \asymp \, \left\{\begin{array}{lll}
(m-n+1)^{\frac12 - \frac 1p} & \qquad & 
\mbox{if}\quad 2 \le p \le \infty\, , \\
1 && \mbox{if}\quad 1 \le p <2\quad \mbox{and}
\quad 1 \le n \le m^{2/p'}\, , \\
m^{1/p'}\, n^{-1/2} && \mbox{if}\quad 1 \le p < 2
\quad \mbox{and}\quad  m^{2/p'} \le n \le m\, . 
\end{array}
\right.
\end{equation}
A  simple monotonicity argument leads to the 
following supplement to  $p=1$.
There exists a constant $c$, independent of $m$ and $n$, such that
\begin{equation}\label{Gluskin*}
d^n (I, \ell_p^m, \ell_2^m) \le c \, n^{-1/2}  
\end{equation}
if $0<  p <1 $ and $1 \le n \le m $.\\ 
The Gelfand widths are  examples of  so-called $s$-numbers, cf.
\cite{Pi87,Pi85} and \cite{CS}. Following Pietsch \cite[2.2.4, p.~80]{Pi87} 
we associate
with  the sequence  of Gelfand widths the following operator ideals.
Let $F$ and $E$ be quasi-Banach spaces and denote by $\cl(F,E)$ the class of
all linear continuous operators $T: \, F \to E$. Then,
for $0<p<\infty$, we put
\[
\cl_{r,\infty}^{(c)} := \Big\{T \in \cl (F,E): \quad \sup_{n \in \Nb}
\, n^{1/r}\, d^n (T) < \infty \Big\}\, .
\] 
Equipped with the quasi-norm 
\[
\lambda_r (T) :=   \sup_{n \in \Nb} \, n^{1/r}\, d^n (T) , 
\]
the set $\cl_{r,\infty}^{(c)}$ becomes a quasi-Banach space.
For such quasi-Banach spaces there always exist a real number 
$\varrho \in (0,1]$ and an equivalent quasi-norm, here denoted by 
$\| \, \cdot \, | \cl_{r,\infty}^{(c)}\|$,
such that
\begin{equation}\label{eq111}
\| \, T_1 + T_2\,  | \cl_{r,\infty}^{(c)}\|^\varrho
\le 
\| \, T_1 \,   | \cl_{r,\infty}^{(c)}\|^\varrho + 
\| \, T_2 \,  | \cl_{r,\infty}^{(c)}\|^\varrho
\end{equation}
holds for all $T_1, T_2 \in \cl_{r,\infty}^{(c)}$. 
\\
To shorten notation we shall use the abbreviation
$I^m_{p,q}$ for the identity $I : \, \ell_p^m \to \ell_q^m$.
It is not complicated to check that (\ref{Gluskin}), (\ref{Gluskin*}) 
imply the following estimates for $\| \, I_{p,2}^m\, 
|\cl^{(c)}_{r,\infty}\| $,
cf. \cite{Li}.

\begin{lemma}\label{Linde}
Let $0 < r < \infty$.\\
{\rm (i)} Let $2 \le p \le \infty$. Then 
\begin{equation}\label{eq110}
\| \, I^m_{p,2} \,| \cl_{r,\infty}^{(c)} \| 
\asymp \, m^{1/r - 1/p + 1/2}
\end{equation}
holds. \\
{\rm (ii)} Let $1  < p < 2$. 
Then 
\begin{equation}\label{eq112}
\| \, I^m_{p,2} \, | \cl_{r,\infty}^{(c)} \| \asymp \, 
\left\{\begin{array}{lll}
m^{1/r - 1/p + 1/2} & \qquad & \mbox{if}\quad 
0 < r \le 2 \, , \\
m^{2/(r p')} && \mbox{if}\quad 2 < r < \infty \, ,
\end{array}\right.
\end{equation}
holds.\\
{\rm (iii)} Let $0  < p \le 1$. 
Then there exists a constant 
$c$ such that
\begin{equation}\label{eq112*}
\| \, I^m_{p,2} \, | \cl_{r,\infty}^{(c)} \| \le c \, 
\left\{\begin{array}{lll}
m^{1/r -  1/2} & \qquad & \mbox{if}\quad 
0 < r \le 2 \, , \\
1 && \mbox{if}\quad 2 < r < \infty \, ,
\end{array}\right.
\end{equation}
holds for all $m \in \Nb$.
\end{lemma}

\noindent
To prove the estimates of the Gelfand numbers from above,  
it turns out to be useful to split the identity $I$ into two parts
$\Id^1, \Id^2$ 
and to treat them independently.
In fact, we shall investigate $\| \, \Id^i\, |\cl^{(c)}_{r_i,\infty}\|$,
$i=1,2$, where $r_1$ and $r_2$ are chosen in different ways. 
For basic properties of the Gelfand numbers we refer to Remark~\ref{r7}
and \cite[2.3]{CS}.
    
\begin{thm}   \label{t7} 
Let $0 < q \le \infty$.\\
{\rm (i)} Let $1 \le p <2$ and  suppose that
$t > d/2$.
Then
\[
d^n (I, b^{s+t}_{p,q},  b^s_{2,2} ) \asymp
n^{-t/d}\,  . 
\]
{\rm (ii)} Let $2< p \le \infty$ and  suppose that $t>0$.
Then 
\[
d^n (I, b^{s+t}_{p,q},  b^s_{2,2} ) \asymp
n^{-t/d }\,  . 
\]
{\rm (iii)} Let $0< p <1$ and  suppose that 
\begin{equation}\label{eq109**}
t > d\, \Big(\frac 1p - \frac 12\Big) \, .
\end{equation}
Then there exist two constants $c_1$ and $c_2$ such that
\[
c_1\, n^{-t/d} \le 
d^n (I, b^{s+t}_{p,q},  b^s_{2,2} ) \le c_2
n^{-t/d  -1 + 1/p}\,  . 
\]
\end{thm} 

\begin{proof}
Without loss of generality we may assume $s=0$.
To see this consider the diagram
\[
\begin{CD}
b^{s+t}_{p,q} @ >I_1>> b^s_{2,2} \\
@V L_{-s} VV @AA L_s A\\
b^{t}_{p,q}@>I_2>> b^0_{2,2}\, ,
\end{CD}
\]
where $L_s$ denotes the isomorphism introduced in
Remark \ref{iso}.
The multiplicativity of the Gelfand numbers implies that 
\[
d^n(I_1,b^{s+t}_{p,q}, b^s_{2,2}) 
\le \| \,  L_{-s} \, \| \, \| \, L_s \, \|\, \, 
d^n (I_2,b^{t}_{p,q}, b^0_{2,2})\, ,
\]
compare with Remark \ref{r7}.
Changing $L_{-s}$ into $L_s$ and vice versa in the diagram above
we end up with
$$
d^n(I_1,b^{s+t}_{p,q}, b^s_{2,2}) = d^n (I_2,b^{t}_{p,q}, b^0_{2,2}).
$$

{\em Step 1.} Estimate from above.
We concentrate on natural numbers $n = 2^{Nd}$
for $N \in \Nb$ (the remaining 
can be treated by the monotonicity of the $d^n$).
Let $\Id_j$ denote the projection given by
\[
\Big(\Id_j \, a\Big)_{m,\lambda} := \left\{ \begin{array}{lll}
a_{j,\lambda} & \qquad &\mbox{if} \quad m =j\, , \\
0 && \mbox{otherwise} \, .
\end{array}\right.
\] 
We split the identity $I$
into a sum $I = \Id^1 + \Id^2$ depending on $N$, where
\[
\Id^1 := \sum_{j=0}^N \Id_j \qquad \mbox{and} \qquad \Id^2 := 
\sum_{j=N+1}^\infty \Id_j\, .
\]
Later on we shall apply the following observation.
Consider the diagram
\[
\begin{CD}
b_{p_,q}^{t} (\nabla) @>\Id_j>> b^0_{2,2} (\nabla)\\
@V P VV @AA Q A\\
\ell_{p}^{|\nabla_j|} @ >I_{p,2}^{|\nabla_j|} >> \ell_{2}^{|\nabla_j|}\, .
\end{CD}
\]
where $P$ and $Q$ are defined as follows. Let 
$a=(a_{\ell,\lambda})_{\ell,\lambda} $. Then
\[
(P(a))_{\lambda} := a_{j,\lambda} \, .
\]
For $b = (b_{\lambda})_{\lambda}$ we define
\[
(Q(b))_{\ell,\lambda} := \left\{ \begin{array}{lll}
a_{j,\lambda} & \qquad & \mbox{if}\quad j=\ell\, , \\
0 && \mbox{otherwise}\, .
\end{array}\right. 
\]
Obviously,
\[
\| \, P \, \| = 2^{-j(t+d(1/2 - 1/p))}\qquad \mbox{and} 
\qquad \| \, Q\, \|= 1\, .
\]
Then property (\ref{mult}) for the Gelfand numbers  yields
\begin{eqnarray}
d^{n} (\Id_j,b^{s+t}_{p,q} , b^s_{2,2})
 & \le &  \| \, P \, \| \, \, \| \, Q\, \|\, \, 
d^{n} (I_{p,2}^{|\nabla_j|}) 
\nonumber
\\
\label{eq113a}
& \le &   2^{-j(t+d(1/2 - 1/p))}\, 
d^{n} (I_{p,2}^{|\nabla_j|} ) \, .
\end{eqnarray}
{\em Substep 1.1.} 
The estimate of 
$d^{n} (\Id^1, b^t_{p,q}, b^0_{2,2} )$, $n=2^{Nd}$.
First we suppose $2 \le p \le \infty$. 
Thanks to (\ref{eq111}), (\ref{eq110}), and 
(\ref{eq113a}) we find
\begin{eqnarray}\label{eq113}
\| \, \Id^1 \, |\cl_{r,\infty}^{(c)} \|^\varrho & \le &  
\sum_{j=0}^N \|\, \Id_j \, |\cl_{r,\infty}^{(c)}\|^\varrho
\nonumber
\\
 & \le &  
\sum_{j=0}^N 2^{-j(t + d(1/2 - 1/p))\varrho}\, 
\|\, I^{|\nabla_j|}_{p,2}  \,|\cl_{r,\infty}^{(c)} \|^\varrho 
\nonumber
\\
& \le &  c_1 \, 
\sum_{j=0}^N \, 2^{-j(t + d(1/2 - 1/p))\varrho}\,  
2^{jd (1/r - 1/p +
  1/2) \varrho}
\nonumber
\\
& \le &  c_2 \, 
2^{ N (d/r - t)\varrho} 
\end{eqnarray}
if $d > t \, r$.
Choosing $r$ small enough, we derive from the definition of 
$\cl_{r,\infty}^{(c)}$ that 
\begin{equation}\label{eq114}
d^n (\Id^1) = d^{2^{Nd}} (\Id^1) \le c_3 \, 2^{- Nt} = c_3 \, n^{-t/d} \, . 
\end{equation}
Now we consider the case  $1 \le p <2 $. 
As above, but using
(\ref{eq112}) instead of (\ref{eq110}), we find
\[
\| \, \Id^1 \, |\cl_{r,\infty}^{(c)} \|    \le   c_2\, 
2^{N (d/r - t)} 
\]
if $1/r > t/d$ and $1/r \ge 2$.
Choosing $r$ small enough,  we obtain
\begin{equation}\label{eq115}
 d^{2^{Nd}} (\Id^1) \le c_4 \, 2^{-Nt}\, . 
\end{equation}
Finally, we investigate the case $0 < p< 1$.
As above,  we obtain 
\begin{equation}\label{eq115*}
d^{2^{Nd}} (\Id^1) \le c_5 \, 2^{-N(t+d-d/p)} 
= c_5 \, n^{-t/d - 1+ 1/p}\, . 
\end{equation}
{\em Substep 1.2.} 
The estimate of 
$d^{n} (\Id^2 , b^t_{p,q}, b^0_{2,2} )$, where  $n=2^{Nd}$. \\
Again we split our considerations into the three 
cases $p\ge 2$ and $1 \le p < 2$ and $0 <p<1$. 
First, let $2 \le p \le \infty$. Using (\ref{eq111}),
(\ref{eq110}), and (\ref{eq113a}),  we find that 
\begin{eqnarray}\label{eq116}
\| \, \Id^2 \, |\cl_{r,\infty}^{(c)} \|^\varrho & \le &  
\sum_{j=N+1}^\infty  
\|\, \Id_j  \, |\cl_{r,\infty}^{(c)} \|^\varrho 
\nonumber
\\
 & \le &  
\sum_{j=N+1}^\infty 2^{-j(t + d(1/2 - 1/p))\varrho}\, 
\|\, I^{|\nabla_j|}_{p,2}  \,|\cl_{r,\infty}^{(c)} \|^\varrho 
\nonumber
\\
& \le &  c_1 \, 
\sum_{j=N+1}^\infty  2^{-j(t + d(1/2 - 1/p))\varrho}\, 
\, 2^{jd (1/r - 1/p + 1/2) \varrho}
\nonumber
\\
& \le &  c_2 \, 2^{ N (d/r - t)\varrho } 
\end{eqnarray}
if $t\, r >d$.
Choosing $r$ large enough ($t >0$ by assumption), 
we derive 
\begin{equation}\label{eq117}
 d^{2^{Nd}} (\Id^2) \le c_3 \, 2^{-Nt}\, . 
\end{equation}
Now we consider  $1 \le p <2 $.  Similarly
\[
\| \, \Id^2 \,  |\cl_{r,\infty}^{(c)} \|\le   c_3\, 
2^{N (d/r - t) } 
\qquad \mbox{if} \quad 
\frac 12 \le \frac 1r < \frac{t}{d} \, .
\]
Since  $t > d/2$,  such a  choice is always possible.
Consequently, 
\begin{equation}\label{eq118}
 d^{2^{Nd}} (\Id^2) \le  c_4 \, 2^{-Nt}\, . 
\end{equation}
Finally, let $0 <p<1$. Then
\begin{equation}\label{eq119}
 d^{2^{Nd}} (\Id^1) \le c_5 \, 2^{-N(t+d-d/p)} \qquad \mbox{if}\qquad 
\frac td + 1- \frac 1p > \frac 1r \ge \frac 12\, . 
\end{equation}
Such a choice is always possible if (\ref{eq109**}) holds.\\
{\em Substep 1.3.}
The additivity of the Gelfand widths yields
\[
d^{2n} (\Id) \le d^n (\Id^1) + d^n (\Id^2) \, .
\] 
In view of this inequality,   the estimate from above of the 
Gelfand widths follows from 
(\ref{eq114})--(\ref{eq119}).\\
{\em Step 2.} Estimate from below. Since $b_n \le c \, d^n$, cf. 
Lemma \ref{bernstein}(i),
we may use Lemma \ref{bernst} here to derive the lower bound in the case
$0< p \le 2$. For $p>2$,   we shall use a 
different  argument. Again we restrict ourselves to a subsequence of the 
natural numbers $n$, where 
\[
\frac{|\nabla_N|}{2} \le n 
<  \frac{|\nabla_N|}{2} +1 \, , \qquad N \in \Nb \, .
\]
Consider the diagram
\[
\begin{CD}
\ell_{p}^{|\nabla_N|}@>I_1>> \ell_{2}^{|\nabla_N|}\\
@V P VV @AA Q A\\
b^{t}_{p,q}(\nabla)@>I_2>> b^0_{2,2}(\nabla)\, ,
\end{CD}
\]
where $I_1$ and $I_2$ denote identities and
this time  $P$ and $Q$ are defined as follows. Let 
$b=(b_\lambda)_{\lambda \in \nabla_N}$. Then
\[
(P(b))_{j,\lambda} := \left\{\begin{array}{lll}
b_\lambda &\qquad & \mbox{if}\quad j=N\, , \\
0 && \mbox{otherwise} \, .
\end{array}\right.
\]
For $a = (a_{j,\lambda})_{j,\lambda}$ we define
\[
(Q(a))_{\lambda} := a_{N,\lambda} \, , \qquad  \lambda \in \nabla_N\, . 
\]
Obviously,
\[
\| \, P \, \| = 2^{N(t+d(1/2 - 1/p))}
\qquad \mbox{and} \qquad \| \, Q\, \|= 1\, .
\]
Then property (\ref{mult}) for the Gelfand numbers  yields that 
\[
d^{n} (I_1,\ell_{p}^{|\nabla_N|} , \ell_{2}^{|\nabla_N|})
 \le  \| \, P \, \| \, \, \| \, Q\, \|\, \, 
d^{n} (I_2,b^{t}_{p,q} (\nabla), b^0_{2,2}(\nabla) ) 
\]
which, in view of Gluskin's estimates (\ref{Gluskin}),  implies that 
\[
c \, 2^{Nd(1/2- 1/p)} 
 \le \,  2^{N(t+ d(1/2-1/p))} \, 
d^{n} (I_2,b^{t}_{p,q} , b^0_{2,2})\, \,   
\]
for some positive $c$ (independent of $N$). This completes the 
estimate from below.
\end{proof}

\begin{rem} 
The use of operator ideals in such a connection
and the associated 
splitting technique applied in Step 1
has some history, cf. {\rm  \cite{C81,Li,KLSS}.} 
Closest to us is {\rm \cite{KLSS}}, where these methods have been used
in connection with entropy numbers. 
\end{rem}


\subsection{Widths of Embeddings of Besov Spaces}   \label{w2} 


Here we do not formulate a general result,  since the restrictions on 
the domains are different for different widths. 


\subsubsection{The Manifold Widths of the Identity } \label{w21}  


\noindent
The main result of this subsection consists in the 
following non-discrete counterpart of Theorem \ref{th3}.

\begin{thm}   \label{th4} 
Let $\Omega $ be a bounded Lipschitz domain.
Let $0 <  p_0,p_1 \le \infty$, $0 <  q_0,q_1 \le \infty$, 
and $s \in \R $. Suppose that {\rm (\ref{einb*})} holds. Then we have 
\begin{equation}\label{eq203*}
e_n^\cont (I,  B^{s+t}_{q_0} (L_{p_0} (\Omega)),   
B^s_{q_1}  (L_{p_1} (\Omega)) )   \asymp 
n^{-t/d}\,  . 
\end{equation}
\end{thm} 

\begin{rem}
Theorem {\rm \ref{th4}}  has several forerunners. We would like to mention
DeVore, Howard, and Micchelli {\rm \cite{DHM89}},
DeVore, Kyriazis, Leviatan, and  Tikhomirov {\rm \cite{DKLT93},} and 
Dung and Thanh {\rm \cite{DD96}.} 
In these papers,  the authors consider the 
quantities $e^\cont_n (I, B^t_{q_0}(L_{p_0}(\Omega)), L_{p_1}(\Omega))$.
Note that from the continuous embeddings 
\[
B^0_{1}(L_p(\Omega)) \hookrightarrow L_p (\Omega)\hookrightarrow 
B^0_{\infty}(L_p (\Omega))\, , \qquad 1 \le p \le \infty\, ,
\] 
we obtain as a direct consequence of Theorem \ref{th4} 
\begin{equation}\label{eq203**}
e_n^\cont (I, B^{t}_{q_0} (L_{p_0} (\Omega)),   L_{p_1} (\Omega))   \asymp 
n^{-t/d}\,  ,
\end{equation}
as long as $1 \le p_1 \le \infty$ and $t>(1/p_0-1/p_1)_+$.
So, Theorem \ref{th4} covers the results obtained before.
However, let us mention that we used the ideas from {\rm \cite{DKLT93}}
for our estimate from above and the ideas from  {\rm \cite{DD96}} to derive 
the estimate from below (here on the level of sequence spaces). 
\end{rem}

\noindent
{\bf Proof of Theorem \ref{th4}}.
Let $\ce $ denote a universal 
bounded linear extension operator  
corresponding to $\Omega$, see Proposition \ref{rychkov} 
in Subsection \ref{lipemin}.
Let $\diam \Omega$ be the diameter of $\Omega$ and let
$x^0 $ be a point in $\R^d$ such that
\[
\Omega \subset \{y: \: |x^0 -y |\le \diam \Omega\} \, .
\] 
Without loss of generality,  we assume that 
\[
\supp \ce f \subset \{y: \: |x^0 -y |\le 2\, \diam \Omega\} \, .
\]
Let $\nabla$ be defined as in (\ref{nabla1}) 
and (\ref{nabla2}) (with $\Omega $ replaced by the ball with radius 
$2 \, \diam \Omega$ and center $x^0$).
Let $R$ denote the restriction operator with respect to $\Omega$.
Let $T$ denote the continuous linear operator that associates to 
$f $ its wavelet series; $T^{-1}$ is the inverse operator. Here
we assume that we can characterize the Besov spaces 
$B^{s+t}_{p_0,q_0} (\R^d)$,  as well as $B^{s}_{p_1,q_1} (\R^d)$, 
in the sense of Proposition \ref{wavelets} in Subsection \ref{app3}.
Then we consider the diagram

\[
 B^{s+t}_{q_0} (L_{p_0}(\Omega))\stackrel{\ce}{\longrightarrow}   
B^{s+t}_{q_0} (L_{p_0}(\R^d))  \stackrel{T}{\longrightarrow} 
b^{s+t}_{p_0,q_0} (\nabla)   
\]
\begin{equation}\label{dia1}
I_1 \hspace{0.5cm} \downarrow \hspace{5cm}
\downarrow \hspace{0.5cm} I_2 
\end{equation}
\[
 B^{s}_{q_1} (L_{p_1}(\Omega)) \stackrel{R}{\longleftarrow}   
B^s_{q_1} (L_{p_1}\R^d)) \stackrel{\, T^{-1}}{\longleftarrow} 
b^{s}_{p_1,q_1} (\nabla) \, .   
\]

\medskip\noindent 
Observe that $I_1 =R \circ T^{-1} \circ I_2 \circ T \circ\ce $.
{}From
(\ref{dia1}) and \eqref{mult} for $e^\cont$, 
we derive that 
\[
e_n^\cont (I_1, B^{s+t}_{q_0} (L_{p_0}(\Omega)),   
B^s_{q_1} (L_{p_1}(\Omega))   \le 
\|\ce \|\, \| \, T\| \, \| T^{-1}\| \, \, 
e_n^\cont (I_2, b^{s+t}_{p_0,q_0}(\nabla) ,   
b^s_{p_1,q_1} (\nabla))   \,  . 
\]
For the converse inequality,  we choose $\nabla^*=(\nabla_j^*)_j$ such that
\[
\supp \psi_{j,\lambda}
\subset \Omega\,, \qquad \lambda \in \nabla_j^*\, , 
\quad j=-1,0,1,\ldots \, , 
\]
and $\inf_j \, 2^{-jd} \, |\nabla_j^*|>0 $.
Then we consider the diagram
\begin{equation}\label{dia2}
\begin{CD}
b^{s+t}_{p_0,q_0}(\nabla^*)@>I_2>> b^{s}_{p_1,q_1}(\nabla^*)
\\
@V T^{-1} VV @AA T A\\
B^{s+t}_{q_0}(L_{p_0}(\Omega))@>I_1>> B^s_{q_1}(L_{p_1}(\Omega))\, ,
\end{CD}
\end{equation}
and conclude that 
\[
e_n^\cont (I_2, b^{s+t}_{p_0,q_0} (\nabla^*), b^s_{p_1,q_1} (\nabla^*))   
\le 
\| \, T\| \, \| T^{-1}\| \, \, 
e_n^\cont (I_1, B^{s+t}_{q_0} (L_{p_0}(\Omega)),   B^s_{q_1} 
(L_{p_1}(\Omega)))   \,  . 
\]
Now Theorem \ref{th3} yields the desired result.
\eproof


\subsubsection{The Widths of Best $m$-Term Approximation of the Identity}   


Let $\Omega$ be a bounded Lipschitz domain in $\Rd$. 
We assume that for any fixed triple $(t,p,q)$ of parameters
the spaces ${B}^{s+t}_{q} (L_p (\Omega))$ and 
$B^{s}_2 (L_2 (\Omega))$ allow a discretization 
by one common wavelet system $\B^*$.
More exactly, we assume that (\ref{eq205})--(\ref{208}) 
are satisfied simultaneously 
for both spaces, cf. Appendix \ref{app8}.
{}From this,  it follows that 
$\B^* \in \B_{C^*}$ for some $1\le C^* < \infty$.

\begin{thm}   \label{th5} 
Let $\Omega$ be as above.
Let $0 < p \le \infty$, $0 < q \le \infty$, $s \in \R $ and  
\[
t> d\, \Big(\frac{1}{p} - \frac 12 \Big)_+
\] 
holds. Then, for any $C\ge C^*$  we have 
\[
e_{n,C}^\non (I, B^{s+t}_{q} (L_p(\Omega)),   
B^s_{2} (L_2(\Omega)))   \asymp 
n^{-t/d}\,  . 
\]
\end{thm} 

\begin{rem}\label{VorLimes} 
i) 
Periodic versions on the $d$-dimensional torus $T^d$ 
may be found in Temlyakov {\rm \cite{Tem00,Tem02}}  with 
$B^s_{2}(L_{2} (\Omega))$ replaced by $L_{p_1} (T^d)$ 
and
$p_1, p, q \ge 1$. 
Furthermore, more general classes of functions are investigated there
(anisotropic Besov spaces, functions of dominating mixed smoothness). 
Finally, let us mention that estimates from below for the quantities 
$$
\inf_{\B \in {\cal O}} 
\sup_{\Vert u \Vert_{B^t_{q_1}(L_{p_1}(T^d))} \le 1 } 
\sigma_n (u, \B)_{L_2 (T^d)} ,
$$
where ${\cal O}$ is the set of all orthonormal bases, have been given by
Kashin $(p_1=q_1=\infty, \ d=1$) and 
Temlyakov {\rm \cite{Tem00,Tem02}}  (general anisotropic case). 
Instead of the manifold widths these authors use entropy numbers. 

ii) We stress that, in this paper, we study the approximation 
in some Hilbertian smoothness space 
$B^s_{2} (L_2(\Omega))$ while most known results from the literature
concern approximation in an  $L_p(\Omega)$-space. 
\end{rem} 

\begin{rem}  \label{Limes} 
We also recall the following limiting case. 
Let
$0 < p < 2$ and $t= d (1/p -  1/2)$. 
Then the embedding $B^{s+t}_{p} (L_p(\Omega)) \hookrightarrow 
B^s_{2} (L_2(\Omega))$ is continuous but not 
compact, cf. Proposition {\rm \ref{compact}.}  Here we have 
\[
\left(\sum_{n=1}^\infty \Big[ n^{t/d} \, 
\sigma_n (u,\B^*)_{B^s_{2} (L_2(\Omega)} \Big]^p \, \frac 1n  \right)^{1/p} 
< \infty \quad \mbox{if and only if} \quad  
u \in  B^{s+t}_{p} (L_p (\Omega)) \, .
\]
A proof can be found in {\rm \cite[Prop.~1]{DDD},} 
but the argument there is mainly based on DeVore and 
Popov {\rm \cite{DP},} see
also {\rm \cite{DJP}.}
\end{rem}

\noindent
{\bf Proof of Theorem \ref{th5}}.
Let $\B^*$ be a wavelet basis as in Appendix \ref{app8}. 
Let $\B$ denote 
the canonical orthonormal basis of $b^0_{2,2} (\nabla)$. 
We equip the Besov space with the equivalent quasi-norm (\ref{208}).
Observe,
\[
\sigma_n (f,\B^*)_{ B^s_{p_1,q_1} (\Omega)} \le c\,
\sigma_n (( \langle f, \widetilde{\psi}_{j,\lambda} 
\rangle )_{j,\lambda}, \B)_{b^s_{p_1,q_1}(\nabla)} 
\, ,
\]
where $c$ is one of the constants in (\ref{eq209}). By means of Theorem 
\ref{th3} and Remark \ref{remark2}(iii), 
this implies the estimate from above. 
The estimate from below follows by combining
Theorem \ref{t1} and Theorem \ref{th4}.
\eproof

\noindent
The simple arguments used in the
proof of Theorem \ref{th5} allow us  to carry over 
Remark~\ref{Limes}
to the sequence space level, see Remark \ref{Limesfall}, 
and Theorem \ref{t5*} to the level of function spaces.

\begin{thm}\label{bmt}
Let $\Omega$ and $\B^*$ be as above.
Let $0 < p_0, p_1, q_0, q_1 \le \infty $, $s\in \R $ and $t>0$
such that {\rm (\ref{einb*})} holds. 
Then we have 
\[
\sup \Big\{  \sigma_n (u,\B^*)_{B^{s}_{q_1} (L_{p_1}(\Omega))} \, :  
\quad \|\, u | B^{s+t}_{q_0} (L_{p_0}(\Omega))
 \| \le 1\Big\} 
 \asymp  n^{-t/d}\,  . 
\]
\end{thm} 

\begin{rem}
\begin{itemize}
\item[i)] 
For earlier results in this direction we refer to 
Kashin {\rm \cite{Ka85}}, Oswald {\rm \cite{Os}}, Donoho {\rm \cite{Do93}}
and DeVore, Petrova and Temlyakov {\rm \cite{DPT03}}.
\item[ii)] 
Not all orthonormal systems are of the same quality, 
see Donoho {\rm \cite{Do93}}. Let us mention the following 
result of DeVore and Temlyakov
{\rm \cite{DT95}.} Let $\B^\# $ denote the 
trigonometric system in $\Rd$. 
By $B^s_q (L_p (\Td))$ we mean the periodic Besov 
spaces defined on the $d$-dimensional torus $\Td$.
Then we put
\[
t(p_0,p_1):= \left\{\begin{array}{lll}
d \big( 1/p_0 - 1/p_1\big)_+ & 
\quad & \mbox{if}\quad 0 < p_0 \le p_1 \le 2 \ \mbox{or}
\ 1 \le  p_1 \le p_0 \le \infty \, , \\
&& \\
d \, \max \big( 1/p_0, 1/2\big) && \mbox{otherwise}\, . 
\end{array}\right.
\]
If $1 \le p_1 \le \infty$, 
$0 < p_0,q_0 \le \infty$, and $t > t(p_0,p_1)$, then
\begin{eqnarray*}
\sup \Big\{ \sigma_n (u,\B^\#)_{L_{p_1} (\Td)} 
: &&\| \, u \, |B^{t}_{q_0}(L_{p_0}(\Td))\| \le 1\Big\} 
\\
 & \asymp & \left\{\begin{array}{lll}
n^{-t/d} & \quad & \mbox{if} \quad p_0 \ge \max (p_1,2)\, , \\
n^{-t/d + 1/p_0 - 1/2} && 
\mbox{if} \quad p_0 \le \max (p_1, 2) =2\, , \\
n^{-t/d + 1/p_0 - 1/p_1} && 
\mbox{if} \quad p_0 \le \max (p_1,2)=p_1 \, .
\end{array}\right.
\end{eqnarray*}
\end{itemize} 
\end{rem}


\subsubsection{The Approximation Numbers  of the Identity}  \label{w34} 


\begin{thm}   \label{th6} 
Let $\Omega$ be a bounded Lipschitz domain.
Let $0 < p \le  \infty$, $0 < q \le \infty$, 
and $s \in \R $. Suppose that
\[
t> d\, \Big(\frac{1}{p} - \frac 12 \Big)_+
\]
holds. 
Then we have 
\[
e_{n}^\lin (I,  B^{s+t}_{q} (L_p (\Omega)),   
B^s_{2}  (L_2 (\Omega)) )    \asymp \left\{
\begin{array}{lll}
n^{-t/d} & \qquad & \mbox{if}\quad 2 \le p \le \infty \, , \\
n^{-t/d + 1/p - 1/2} & 
\qquad & \mbox{if}\quad 0 < p <2 \,  . 
\end{array}\right.
\]
\end{thm} 

\begin{proof}
The statement is a consequence of Theorem \ref{th3}(ii), 
Proposition \ref{rychkov}, (\ref{eq203d}) and (\ref{eq204}). 
\end{proof}

\begin{rem} \label{partial}
\begin{itemize}
\item[(i)]
The proof is constructive.
An order-optimal linear approximation is obtained by taking an appropriate
partial sum of the wavelet series of ${\mathcal E}f$, where 
${\mathcal E}$ is the linear universal extension operator from
Proposition {\rm \ref{rychkov}}, cf. 
Remark {\rm \ref{opt}} for the discrete case.
\item[(ii)]
This result is well-known. It can be derived from  
{\rm \cite{T02}} and {\rm \cite[3.3.2]{ET96}.}
There and in {\rm \cite{Cae}}  information 
can also  be found about what is known 
for the general situation, i.e., in which  $B^s_{2} (L_2(\Omega))$
is replaced by 
$B^s_{q_1} (L_{p_1}(\Omega))$.
However, let us mention that there are many references 
which had dealt with this problem before;
we refer to {\rm \cite[Thm.~1.4.2]{Tem93}} 
and {\rm \cite[Thm.~9,~p.193]{Ti90}} and the comments given there.
\end{itemize}
\end{rem}


\subsubsection{The Gelfand  Widths of the Identity}   



\begin{thm}   \label{t10} 
Let $\Omega \subset \Rd$ be a 
bounded Lipschitz domain  and let $0 < q \le \infty$.\\
{\rm (i)} Let $1 \le p <2$ and  suppose that 
$t > d/2$.
Then
\[
d^n (I, B^{s+t}_{q}(L_p(\Omega)),  B^s_{2}(L_2(\Omega)) ) \asymp
n^{-t/d}\,  . 
\]
{\rm (ii)} Let $2< p \le \infty$ and  suppose that  $t>0$.
Then 
\[
d^n (I, B^{s+t}_{q}(L_p(\Omega)),  B^s_{2}(L_2(\Omega)) ) \asymp
n^{- t/d }\,  . 
\]
{\rm (iii)} Let $0< p <1$ and  suppose that 
\[
t > d\, \Big(\frac 1p - \frac 12\Big) \, .
\]
Then there exists two constants $c_1$ and $c_2$ such that
\[
c_1\, n^{-t/d} \le 
d^n (I, B^{s+t}_{q}(L_p(\Omega)),  B^s_{2}(L_2(\Omega)) ) \le c_2
n^{-t/d  -1 + 1/p}\,  . 
\]
\end{thm} 

\begin{proof}
Consider the diagram
\[
\begin{CD}
B^{s+t}_{q_0}(L_{p_0}(\Omega))@>I_1>> B^s_{2}(L_2(\Omega))
\\
@V T VV @AA T^{-1} A\\
b^{s+t}_{p_0,q_0}(\nabla)@>I_2>> b^s_{2,2}(\nabla)\, ,
\end{CD}
\]

\medskip\noindent 
where $T$ and $T^{-1}$ are defined as in the proof of Theorem 
\ref{th4}. Since 
$I_1= T^{-1}\circ I_2 \circ  T$, 
it is enough to combine property \eqref{mult} for the Gelfand numbers and 
Theorem \ref{t7} to derive the estimates from above. For the estimates 
from below,  one uses the diagram
\[
\begin{CD}
b^{s+t}_{p_0,q_0}(\nabla^*)@>I_1>> b^s_{2,2}(\nabla^*)
\\
@V T VV @AA T^{-1} A\\
B^{s+t}_{q_0}(L_{p_0}(\Omega))@>I_2>> B^s_{2}(L_2(\Omega))\, ,
\end{CD}
\]

\medskip\noindent 
where $\nabla^*$ is defined  as in proof of Theorem 
\ref{th4}. This completes the proof.
\end{proof}

\begin{rem}
Partial results concerning Gelfand numbers of embedding operators
may be found in the monographs Pinkus 
{\rm \cite[Chapt.~VII, Thm.~1.1]{Pi85},} 
Tikhomirov {\rm \cite[Thm.~39,~p.~206]{Ti90},} and
Triebel {\rm \cite[4.10.2]{T78}.} Let $T$ be a compact operator in 
$\cl (F,E)$, where $F,E$ are arbitrary Banach spaces and let
$d_n (T,F,E)$ denote the Kolmogorov numbers. Then 
\[
d^n (T') = d_n (T)\, , \qquad n \in \Nb\, , 
\]
holds, cf. {\rm \cite[Prop.~2.5.6]{CS}}
or {\rm \cite{Pi74}}. 
For Kolmogorov numbers the asymptotic behaviour is also known in 
certain situations, cf.  {\rm \cite[Chapt.~VII, Thm.~1.1]{Pi85},} 
{\rm \cite[Thm.~10,~p.~193]{Ti90},} 
{\rm \cite[4.10.2]{T78},} and {\rm \cite{Tem93}.}
\end{rem}


\subsection{Proofs of Theorems \ref{t2}, \ref{t3}, 
and \ref{t5}}   \label{w3} 


\subsubsection{Proof of Theorem \ref{t2}}      \label{w31} 


For  $s>0$ we have $H^{-s} (\Omega) = 
B^{-s}_2 (L_2(\Omega))$. Hence, Theorem \ref{th6} yields
\[
e_{n}^\lin (I,  B^{-s+t}_{q} (L_p (\Omega)), H^{-s} (\Omega) ) 
\asymp \left\{
\begin{array}{lll}
n^{-t/d} & \qquad & \mbox{if}\quad 0 < p \le 2\, , \\
n^{-t/d + 1/p -1/2} & 
\qquad & \mbox{if}\quad 2 < p \le \infty\,  . 
\end{array}\right.
\]
Since $S: H^{-s} (\Omega) \to H^s_0 (\Omega) $ is an isomorphism, 
we obtain the desired result
from property \eqref{mult} for 
the approximation numbers. 


\subsubsection{Proof of Theorem \ref{t3} }   \label{w32} 


Since of $H^{-s} (\Omega) = B^{-s}_2 (L_2(\Omega))$, 
Theorem \ref{th5} yields that 
\[
e_{n,C}^\non (I,  B^{-s+t}_{q} (L_p (\Omega)),   H^{-s}(\Omega) )    
\asymp n^{-t/d} 
\]
Since $S: H^{-s} (\Omega) \to H^s_0 (\Omega) $ is an isomorphism,  
Lemma \ref{basic}(ii) implies the desired result.  

\subsubsection{Proof of Theorem \ref{t5}}  \label{prooft5} 

All what we need from the wavelet basis
is the following estimate for the best $n$-term approximation
in the $H^1$-norm:
\begin{equation} \label{H1scale}
\|\, u-S_n(f)\, \|_{H^1(\Omega)} \leq c\,  
\|\, u\, |{B^{t+1}_{\tau}(L_{\tau}(\Omega))}\| \, n^{-t/2}, \quad 
\hbox{where}  \quad 
\frac{1}{\tau}= \frac{t}{2} + \frac{1}{2},
\end{equation}
see, e.g., \cite{DDD} (however we could instead use Theorem \ref{bmt}).
We therefore have to estimate
the Besov norm $B^{\alpha}_{\tau}(L_{\tau}(\Omega))$. 
Since $1< p \le 2$,  the
embedding 
\linebreak $B_p^{k-1}(L_p(\Omega))\hookrightarrow W^{k-1}_p(\Omega)$ 
holds, cf. e.g. \cite[2.3.2,~2.5.6]{T83}.
Hence our right--hand side $f$ is contained
in the Sobolev space $W^{k-1}_p(\Omega)$.
Therefore we may employ  the fact that $u$ can be decomposed 
into a regular part $u_R$ and a singular part $u_S$, i.e., 
$u=u_R + u_S,$ where $u_R \in W^{k +1}_p(L_p(\Omega))$ 
and $u_S$  only depends on the shape of 
the domain and can be computed explicitly, cf. 
Grisvard \cite[Thm.~2.4.3]{Gr1}. We introduce
polar coordinates $(r_l, \theta_l)$ in the 
vicinity of each vertex $\Upsilon_l$ and introduce the functions
\[
{\cal S}_{l,m}(r_l,\theta_l) : = \left\{ \begin{array}{lll} 
&& \hspace{-0.5cm}\zeta_l(r_l)r_l^{\lambda_{l,m}}
\sin(m\pi \theta_l /\omega_l)  \qquad  \mbox{if}\quad  
\lambda_{l.m}:=m\pi/\omega_l \neq \mbox{integer}\, , \\
&& \\ 
&& \hspace{-0.5cm}
\zeta_l(r_l)r_l^{\lambda_{l,m}}[\log r_l \sin(m\pi \theta_l/\omega_l)+
\theta_l\cos(m\pi \theta_l/\omega_l)] \qquad \mbox{otherwise}\, .
\end{array}\right.
\]
Here $\zeta_1, \ldots , \zeta_N$ denote  
suitable $C^{\infty}$ truncation functions and $m$ is a natural number.
  Then for $f \in W^{k-1}_p(\Omega)$, one has
  \begin{equation} \label{singpart}
u_S= \sum_{l=1}^N \sum_{0<\lambda_{l,m}< k+1-2/p}c_{l,m} \, 
{\cal S}_{l.m}\, ,
  \end{equation}
  provided that  no $\lambda_{l,m}$ is equal to $k+1-2/p$. This means that
the finite number of singularity functions that is 
needed depends on the scale of  spaces
  we are interested in, i.e., on the smoothness parameter $k$. 
  According to  (\ref{H1scale}), we have to estimate the
  Besov regularity of both, $u_S$ and $u_R$, in the specific scale
  $$
  B^{t+1}_{\tau}(L_{\tau}(\Omega)) , \quad \hbox{where} 
  \quad \frac{1}{\tau}= \frac{t}{2} + \frac{1}{2}\, .
$$
Since $u_R \in W^{k+1}_p(\Omega)$, the boundedness of $\Omega$ implies
the embedding
\[ 
W^{k+1}_p(\Omega) \hookrightarrow {B^{k+1-\delta}_{q}(L_{q}(\Omega))},  
\quad \hbox{with} \quad 
\delta > 0\, , \quad 0 < q \le p\, , \quad
k+1 > 2\, \Big(\frac 1q - \frac 12\Big)\, .
\]
Hence
\begin{equation} \label{regreg}
u_R \in {B^{k+1-\delta}_{\tau}(L_{\tau}(\Omega))},
\quad \hbox{with} \quad 
\frac{1}{\tau}= \frac{(k-\delta)}{2} 
+ \frac{1}{2}\qquad \mbox{for arbitrarily small} \quad \delta >0\, .
\end{equation}
Moreover, it has been shown in \cite{Da1} (see also 
Remark \ref{singularity}) that the 
functions  ${\cal S}_{l,m}$ defined above satisfy 
\begin{equation} \label{singreg} 
{\cal S}_{l,m}(r_l,\theta_l)\in 
{B^{1/2 + 2/q}_{q}(L_{q}(\Omega))},  
\qquad \mbox{for all} \quad 0 < q <\infty\, .
\end{equation}
By combining (\ref{regreg}) and (\ref{singreg}) we see that 
\[
u \in 
{B^{k+1-\delta}_{\tau}(L_{\tau}(\Omega))}, \quad 
\hbox{where} \quad 
\frac{1}{\tau}= \frac{(k-\delta)}{2}+ \frac{1}{2}
\qquad \mbox{for arbitrarily small} \quad \delta >0.
\]
To derive an estimate uniformly with respect to 
the unit ball in $B^{k-1}_p (L_p(\Omega))$
we argue as follows. We put
\[
{\mathcal N} := \span \Big\{ {\cal S}_{l,m}(r_l,\theta_l): 
\quad 0 < \lambda_{m,l}< k+1-2/p\, , 
\: l=1,\,  \ldots \, , N\Big\}\, .
\]
Let $\gamma_l$ be the trace operator with respect to the segment $\Gamma_l$.
Grisvard has shown that
$\Delta$ maps 
\[
H := \Big\{ u \in W^{k+1}_p(\Omega): 
\quad \gamma_l u =0\, , \, l = 1, \ldots \, ,N\Big\} \: + \: {\mathcal N} 
\]
onto $W^{k-1}_p(\Omega)$, cf. \cite[Thm.~5.1.3.5]{Gr3}. 
This mapping is also injective, see\linebreak
 \cite[Lemma~4.4.3.1, Rem.~5.1.3.6]{Gr3}.
We equip the space $H$ with the norm
\[
\| \, u\,  \|_H := \| \, u_R + u_S\,  \|_H =
\| \, u_R\,  \|_{W^{k+1}_p(\Omega)} 
+ \sum_{l=1}^N \sum_{0<\lambda_{l,m}< k+1-2/p} \, |c_{l,m}|\, , 
\] 
see (\ref{singpart}). Then $H$ becomes a Banach space. Furthermore, 
$\Delta : H \to W^{k-1}_p (\Omega) $ is continuous.
Banach's continuous inverse theorem implies  
that the solution operator is continuous,  considered 
as a mapping  from $W^{k-1}_p(\Omega)$ onto $H$.
Finally, observe that 
\[
\| \, u_R + u_S \, \|_{B^{k+1-\delta}_{\tau}(L_{\tau}(\Omega))}  \le
C\, \Big(
\| \, u_R \, \|_{W^{k+1}_p(\Omega)} 
+ \sum_{l=1}^N \sum_{0<\lambda_{l,m}< k+1-2/p} \, |c_{l,m}|\Big) 
\]
with some constant $C$ independent of $u$. 
\eproof


\section{Appendix -- Besov spaces}


Here we collect some properties of Besov spaces 
that have been used in the text before. 
Detailed references will be given. 
For general information on Besov spaces,  we refer 
to the monographs \cite{Me,Ni,Pe,RS96,T83,T92}.


\subsection{Besov Spaces on $\Rd$ and Differences}\label{app1}


Nowadays Besov spaces are widely used in several branches of mathematics.
Probably the most common way to introduce 
these classes makes use of differences.
For $M\in \Nb$, $h\in \Rd$, and $f: \Rd \to \Cd$ we
define
\[
\Delta_h^M f (x) :=
\sum\limits_{j=0}^M {M \choose j} (-1)^{M-j}\, f(x+jh).
\]
Let $0 < p \le \infty$.
The corresponding modulus of smoothness is then given by
\[
\omega^M (t,f)_p := \sup_{|h|<t} \, \| \, 
\Delta_h^M f \,  \|_{L_p (\Rd)}\, , \qquad t>0 \, .
\]
One approach to introduce Besov spaces is the following.

\begin{definition}
Let $s > 0$ and $0 < p,q \le \infty$. 
Let $M$ be a natural number satisfying $M>s$.
Then 
$\Lambda^s_q(L_p (\Rd))$ is the collection of all 
functions  $f \in L_p (\Rd) $ such that
\[
| \, f \, |_{\Lambda^s_q(L_p (\Rd))} := \bigg( \int_0^\infty \Big[
t^{-s} \, \omega^M (t,f)_p\Big]^q \frac{dt}{t}\bigg)^{1/q} <\infty
\]
if $q< \infty$ and
\[
| \, f \, |_{\Lambda^s_\infty(L_p (\Rd))} :=  \sup_{t>0}\, 
t^{-s} \, \omega^M (t,f)_p < \infty
\]
if $q= \infty$. 
These classes are equipped with a quasi-norm by taking
\[
\| \, f \, \|_{\Lambda^s_q(L_p (\Rd))} 
:= \| \, f\, \|_{L_p(\Rd)} +  | \, f \, |_{\Lambda^s_q(L_p (\Rd))}\, .
\]
\end{definition}

\begin{rem}
It turns out that these classes do not depend 
on $M$, cf. {\rm \cite{DVS}.}
\end{rem}

\begin{rem}\label{singularity}
Let $\varrho \in C^\infty_0 (\Rd)$ be a function 
such that $\varrho (0) \neq 0$. 
By means of the above definition it is not 
complicated to show that a function
\[
f_\alpha (x):= |x|^\alpha \, \varrho (x)\, , 
\qquad x \in \R^d\, , \quad \alpha >0\, , 
\]
belongs to $\Lambda^{\alpha + d/p}_\infty(L_p (\Rd))$ 
and that this is best the possible 
(if $\alpha$ is not an even natural number), cf. 
{\rm \cite[2.3.1]{RS96}} for details.
A minor modification shows that 
\[
f_{\alpha,\beta} (x):= |x|^\alpha \, (\log |x|)^\beta\,  
\varrho (x)\, , \qquad x \in \R^d\, , \quad \alpha, \, 
\beta  >0\, , 
\]
belongs to $\Lambda^{\alpha + d/p-\varepsilon}_\infty(L_p (\Rd))$ 
for all $\varepsilon$, $ 0 < \varepsilon < \alpha + d/p$.
\end{rem}


\subsection{Besov Spaces on $\Rd$ and Littlewood-Paley Characterizations}
\label{app2}


Since we are using also spaces with negative smoothness 
$s<0$ and/or $p,q<1$  we shall give a further  
definition, which relies on  Fourier analysis.
We use it here for introductory purposes.
This approach makes use of smooth dyadic decompositions of unity.
Let $\phi \in C^\infty_0 (\Rd)$ be a function 
such that $\phi (x)=1$ if $|x|\le 1$ and 
$\phi (x)=0$ if $|x|\ge 2$. Then we put
\begin{equation}
\phi_0 (x):= \phi(x), \qquad \phi_j (x) 
:= \phi (2^{-j}x) - \phi (2^{-j+1}x)\, , \quad j \in \Nb\, . 
\end{equation}
It follows
\[
\sum_{j=0}^\infty  \phi_j (x) = 1\, , \qquad x \in \Rd\, , 
\]
and \[
\supp \, \phi_j \subset \Big\{x \in \Rd: \quad 2^{j-2} 
\le |x| \le 2^{j+1}\Big\}\, , 
\qquad j=1,2, \ldots \, .
\] 
Let $\cf$ and $\cfi$ denote the Fourier transform and its inverse,  
both defined on $\S'(\Rd)$.
For $f \in \S'(\Rd)$ we consider the sequence 
$\cfi [\phi_j (\xi) \, \cf f (\xi)] (x)$, $j \in \Nb_0$, 
of entire analytic functions.
By means of these functions,  we define the Besov classes.

\begin{definition}
Let $s \in \R$ and $0 < p,q \le \infty$. Then 
$B^s_q(L_p (\Rd))$ is the collection of all 
tempered distributions $f$ such that
\[
\| \, f \, | B^s_q(L_p (\Rd))\|= \bigg( \sum_{j=0}^\infty 
2^{sjq}\, 
\|  \cfi [\phi_j (\xi) \, \cf f (\xi)] 
(\, \cdot \, )\, |L_p (\Rd)\|^q\bigg)^{1/q} <\infty
\]
if $q< \infty$ and
\[
\| \, f | B^s_\infty(L_p (\Rd))\|=  \sup_{j=0,1, \ldots}\, 2^{sj}\,  \|  
\cfi [\phi_j (\xi) \, \cf f (\xi)] (\, \cdot \, )\, |L_p (\Rd)\| <\infty
\]
if $q= \infty$. 
\end{definition}

\begin{rem}\label{confusion}
\begin{itemize}
\item[i)]
If no confusion is possible we drop $\Rd$ in notations.
\item[ii)]
These classes are quasi-Banach spaces. 
They do not depend on the chosen function $\phi$
(up to equivalent quasi-norms). If $t= \min (1,p,q)$, then
\[
\| \, f+g \,| B^s_q(L_p)\|^t \le \| \, f \, | B^s_q(L_p)\|^t + 
\| \, g \,  | B^s_q(L_p)\|^t 
\]
holds for all $f,g \in B^s_q (L_p)$.
\end{itemize}
\end{rem}

\begin{proposition}{\rm \cite[2.5.12]{T83}.}
Let $0 < p,q\le \infty$ and 
$s> d \, \max (0, 1/p -1)$. Then we
have coincidence of $\Lambda^s_q(L_p)$ 
and $B^s_{q} (L_p)$ in the sense of equivalent quasi-norms.
\end{proposition}

\begin{rem} \label{multiplier}
\begin{itemize}
\item[i)]
For $s \le  d \, \max (0, 1/p -1)$ 
we have $\Lambda^s_q(L_p) \neq B^s_{q} (L_p)$. E.g., 
the Dirac distribution $\delta$ belongs 
to $B^{d(1/p-1)}_\infty (L_p)$, cf. {\rm \cite[2.3.1]{RS96}.}
\item[ii)]
Smooth cut-off functions are pointwise multipliers for all Besov spaces. 
More exactly, let 
$\psi \in \D$. Then the product $\psi \, f$ belongs to 
$B^s_q(L_p)$ for any $f \in B^s_q(L_p) $ and there exists a constant $c$ 
such that
\[
\| \, \psi \, f \, |B^s_q(L_p)\| \le c \| \, f \, | B^s_q(L_p)\| 
\] 
holds, see e.g. {\rm \cite[2.8]{T83}, \cite[4.7]{RS96}.}
\end{itemize}
\end{rem}


\subsection{Wavelet Characterizations}
\label{app3}


For the  construction of  biorthogonal wavelet  bases  as considered below, 
we refer to the recent monograph of 
Cohen \cite[Chapt.~2]{C03}. 
Let $\varphi$ be a compactly supported 
scaling function of sufficiently high regularity and let 
$\psi_i$, where  $i=1, \ldots \, 2^d-1$,  be the corresponding wavelets.
More exactly, we suppose for some $N>0$ and $r\in \Nb$
\begin{eqnarray*}
&& \supp \, \varphi\, , \: \supp \, \psi_i 
\quad \subset \quad  [-N,N]^d\, , 
\qquad i=1, \ldots , 2^d-1 \, , \\
&& \varphi,  \psi_i \in C^r(\Rd)\, , \qquad i=1, \ldots , 2^d-1 \, ,   \\
&&
\int x^\alpha \, \psi_i (x)\, dx =0 \qquad \mbox{for all} \quad
|\alpha|\le r\, , \qquad i=1, \ldots , 2^d-1 \, ,  
\end{eqnarray*}
and
\[
\varphi (x - k), \, 2^{jd/2}\, \psi_i (2^jx-k)\, , \qquad j\in \Nb_0\, ,
\quad k \in \Zd \, , 
\]
is a Riesz basis in $L_2(\Rd)$.
We shall use the standard abbreviations
\[
\psi_{i,j,k} (x) = 2^{jd/2}\, \psi_i (2^jx-k) \qquad \mbox{and} \qquad
\varphi_k (x) = \varphi (x-k) \, .
\]
Further, 
the dual Riesz basis should fulfill  the same requirements, i.e., 
there exist functions $\widetilde{\varphi} $ and 
$\widetilde{\psi}_i$, $ i=1, \ldots , 2^d-1$,
such that 
\begin{eqnarray*}
\langle \widetilde{\varphi}_k, \psi_{i,j,k} \rangle   
& = & \langle \widetilde{\psi}_{i,j,k}, 
\varphi_k \rangle  = 0\, , \\
\langle \widetilde{\varphi}_k, \varphi_{\ell} \rangle  
& = & \delta_{k,\ell} \qquad
(\mbox{Kronecker symbol})\, , \\
\langle \widetilde{\psi}_{i,j,k}, 
\psi_{u,v,\ell} \rangle   & = & \delta_{i,u}\,
 \delta_{j,v}\, \delta_{k,\ell}\, , \\
\supp \, \widetilde{\varphi} \, , &&  \supp \, \widetilde{\psi}_i \quad 
\subset \quad  [-N,N]^d\, , 
\qquad i=1, \ldots , 2^d-1 \, , \\
\widetilde{\varphi},   \widetilde{\psi}_i & \in 
&  C^r(\Rd)\, , \qquad i=1, \ldots , 2^d-1 \, ,   \\
\int x^\alpha \, 
\widetilde{\psi}_i (x)\, dx & = & 0 \qquad \mbox{for all} \quad
|\alpha|\le r\, , \qquad i=1, \ldots , 2^d-1 \, .  
\end{eqnarray*}  
For $f \in \S'(\Rd)$ we put
\begin{equation}\label{eq100}
\langle f,\psi_{i,j,k} \rangle  
= f(\overline{\psi_{i,j,k}}) \qquad \mbox{and}\qquad
 \langle f,\varphi_k \rangle  = f(\overline{\varphi_k}) \, ,
\end{equation}
whenever this makes sense.

\begin{proposition}\label{wavelets}
Let $s \in \R$ and $0 < p,q \le \infty$. 
Suppose 
\begin{equation}\label{eq101}
r > \max \Big (s, \frac{2d}{p} + \frac d2 - s\Big)\, .
\end{equation}
Then $B^s_q(L_p)$ is the collection of all 
tempered distributions $f$ such that
$f$ is representable as
\[
f = \sum_{k \in \Zd} a_{k}\, \varphi_k
 + \sum_{i=1}^{2^d-1} \, \sum_{j=0}^\infty 
\sum_{k \in \Zd} a_{i,j,k}\, \psi_{i,j,k} 
\qquad (\mbox{convergence in} \quad \S')
\]
with
\[
\| \, f\, |B^s_q(L_p)\|^* := \Big(\sum_{k \in \Zd} |a_{k}|^p\Big)^{1/p} +
 \bigg(\sum_{i=1}^{2^d-1} \, \sum_{j=0}^\infty 
2^{j(s + d(1/2-1/p))q} 
\Big(\sum_{k \in \Zd} |a_{i,j,k}|^p\Big)^{q/p}\bigg)^{1/q} < \infty\, ,
\]
if $q < \infty$ and
\[
\| \, f\, |B^s_\infty (L_p)\|^* := 
\Big(\sum_{k \in \Zd} |a_{k}|^p\Big)^{1/p} +
  \sup_{i=1, \ldots \, , 2^d-1} \, \sup_{j=0, \ldots \, } \, 
2^{j(s + d(1/2-1/p))} 
\Big(\sum_{k \in \Zd} |a_{i,j,k}|^p\Big)^{1/p} < \infty\, .
\]
The representation is unique and
\[
a_{i,j,k} =  \langle f,
\widetilde{\psi}\begin{normalsize}\end{normalsize}_{i,j,k}\rangle  
\qquad \mbox{and} \qquad 
a_k = \langle f,\widetilde{\varphi}_{k} \rangle  
\] 
hold. Further $I : f \mapsto \{\langle f,
\widetilde{\varphi}_{k} \rangle ,\,  \langle f,
\widetilde{\psi}_{i,j,k}\rangle  \} $
is an isomorphic map of $B^s_q(L_p (\Rd))$ onto 
the sequence space equipped with the quasi-norm 
$\| \, \cdot \, |B^s_q(L_p)\|^*$, i.e.,  
$\| \, \cdot \, |B^s_q(L_p)\|^*$ may serve as an 
equivalent quasi-norm on $ B^s_q(L_p )$.
\end{proposition}

\begin{rem} \label{remark31}
\begin{itemize}
\item[i)] 
The restriction {\rm (\ref{eq101})} 
guarantees that {\rm (\ref{eq100})} makes sense
for all $f \in B^s_q(L_p)$.
\item[ii)] 
It is immediate from this proposition
that the functions $\varphi_k, \psi_{i,j,k}, k \in \Zd, 1 \le i \le 2^d-1,
j\in \Nb_0$ form a basis for $B^s_q(L_p)$ 
if $\max(p,q)< \infty$. By the same
reasoning the functions 
\[\varphi_k, \quad 2^{-js}\, 
\psi_{i,j,k}, \qquad  k \in \Zd, \quad 1 \le i \le 2^d-1,\quad
j\in \Nb_0\, , 
\]
form a Riesz basis for $B^s_2(L_2)$.  
\item[iii)] 
If the wavelet basis is orthonormal (in $L_2$), 
then this proposition is proved in Triebel {\rm \cite{T03}.} 
But the comments made in Subsection {\rm 3.4} of the quoted paper make clear
that this extends to the situation considered in 
Proposition {\rm \ref{wavelets}.}
A different proof, but restricted to 
$s> d(1/p - 1)_+$, is given in {\rm \cite[Thm.~3.7.7]{C03}.}
However, there are many forerunners with some 
restrictions concerning $s,p$ and $q$. 
We refer to   {\rm \cite{Bou}} and {\rm \cite{Me}.} 
\end{itemize}
\end{rem} 


\subsection{Besov Spaces on Domains -- the 
Approach via Restrictions}\label{app4}


There are at least two different approaches to define  function spaces 
on domains. One approach  uses 
restrictions to $\Omega$ of functions defined  on $\Rd$.
So, all calculations are done on $\Rd$.
The other approach introduces theses spaces by 
means of local quantities defined only 
in $\Omega$. For numerical purposes the second 
approach is more promising whereas for analytic investigations 
the first one looks more elegant. 
Here we discuss both, since both were used.\\ 
Let $\Omega \subset \Rd$ be an  bounded open nonempty set. Then we define
$B^s_q (L_p(\Omega))$ to be the collection of 
all distributions $f \in \D'(\Omega)$ such that there
exists a tempered distribution $g \in B^s_q (L_p (\Rd))$ satisfying 
\[
f (\varphi) = g(\varphi) \qquad \mbox{for all} 
\quad \varphi \in \D(\Omega) \, ,
\]
i.e. $g|_\Omega = f$ in $\D'(\Omega)$.
We put
\[
\| \, f\, |B^s_q (L_p (\Omega))\|:= 
\inf \, \| \, g \, |B^s_q (L_p (\Rd))\|\, , 
\] 
where the infimum is taken with respect to all distributions $g$ as above.
\\
Let
$\diam \, \Omega$ be the diameter of the 
set $\Omega$ and let $x^0$ be a point with the property
\[
\Omega \subset \Big\{y: \: |x^0-y|\le \diam \Omega\Big\} \, .
\]
Such a point we shall call a \emph{center}  of $\Omega$.
Since smooth cut-off functions are pointwise 
multipliers, cf. Remark \ref{multiplier}, we can associate with any $f \in 
B^s_q (L_p (\Omega))$ a tempered distribution $g \in B^s_q(L_p)$ such that
$g|_\Omega = f$ in $\D'(\Omega)$,
\begin{eqnarray}\label{eq102}
C \, \| \, g\, |B^s_q (L_p )\|  
& \le &    \| \, f \, |B^s_q (L_p (\Omega))\| 
\le  \| \, g\, |B^s_q (L_p )\|  
\\
\label{eq103}
\supp \,  g  & \subset &  \{x \in \Rd : \quad |x-x^0| 
\le 2 \, \diam \, \Omega\}\, .
\end{eqnarray}
Here $0 < C <1$ does not depend on $f$ (but on $\Omega, s, p, q$).\\
Now we turn to decompositions by means of wavelets. 
We use the notation from the preceeding subsection.
Define 
\begin{equation}\label{lambda}
\Lambda_j := \Big\{ k \in \Zd: \quad |k_i-x^0_i|\le 2^j \, \diam \, 
\Omega + N\, , \: i=1, \ldots \, , d \Big\}\, , \qquad j=0,1, \ldots \, .
\end{equation}
Then given $f$ and taking $g$ as above,  we find that 
\begin{equation}\label{eq201}
g = \sum_{k \in \Lambda_0} \langle g, 
\widetilde{\varphi}_k \rangle \, \varphi_k
 + \sum_{i=1}^{2^d-1} \, \sum_{j=0}^\infty 
\sum_{k \in \Lambda_j} \langle g, 
\widetilde{\psi}_{i,j,k} \rangle\, \psi_{i,j,k} 
\qquad (\mbox{convergence in} \quad \S')
\end{equation}
and 
\begin{eqnarray}\label{eq202}
\| \, g\, |B^s_q (L_p)\| & \asymp &  \Big(\sum_{k \in \Lambda_0} 
| \langle g, \widetilde{\varphi}_k \rangle |^p\Big)^{1/p}  +
\\
&&
\bigg(\sum_{i=1}^{2^d-1} \, \sum_{j=0}^\infty \, 
2^{jq(s + d(\frac{1}{2}-\frac{1}{p}))} 
\Big(\sum_{k \in \Lambda_j} | \langle g, 
\widetilde{\psi}_{i,j,k} \rangle |^p\Big)^{q/p}\bigg)^{1/q} < \infty\, .
\nonumber
\end{eqnarray}
The following more handy notation is also used.
We put
\begin{eqnarray}\label{nabla1}
\nabla_{-1} & := & \Lambda_0 \\
\label{nabla2}
\nabla_j & := & \Big\{
(i,k) : \quad 1 \le i \le 2^d-1\, , \quad k \in \Lambda_j
\Big\}\, , \quad j=0,1,\ldots \, , 
\end{eqnarray}
$\psi_{j,\lambda}  := \psi_{i,j,k}$, if $\lambda = (i,k) \in \nabla_j$,
$j \in \Nb_0$, and 
$\psi_{j,\lambda} 
:= \varphi_k $ if $\lambda =k \in \nabla_{-1}$. For the dual basis, 
(\ref{eq201}) and (\ref{eq202}) read as
\begin{equation}\label{eq203d}
g = \sum_{j=-1}^\infty 
\sum_{\lambda \in \nabla_j} \langle g, 
\widetilde{\psi}_{j,\lambda} \rangle\, \psi_{j,\lambda} 
\qquad (\mbox{convergence in} \quad \S')
\end{equation}
and 
\begin{equation}\label{eq204}
\| \, g\, |B^s_q (L_p)\|  \asymp
  \bigg( \sum_{j=-1}^\infty \, 
2^{jq(s + d(\frac{1}{2}-\frac{1}{p}))} 
\Big(\sum_{\lambda \in \nabla_j} | \langle g, 
\widetilde{\psi}_{j,\lambda} \rangle |^p\Big)^{q/p}\bigg)^{1/q} < \infty\, .
\end{equation}

\subsection{Lipschitz Domains, Embeddings, and Interpolation} 
\label{lipemin}

We call a domain $\Omega$ a {\em special 
Lipschitz domain} (see Stein \cite{St}),
if $\Omega$ is an open set in $\Rd$ and
if there exists a function $\omega: \, \R^{d-1} \to \R$ such that
\[
\Omega = \Big\{(x',x_d) \in \Rd : \: x_d > \omega (x') \Big\}
\] 
and 
\[
|\, \omega (x') - \omega (y')\, | \le C \, |x' - y'| \qquad \mbox{for all}
\quad x', \, y' \in \R^{d-1} \, ,
\]
and some constant $C>0$.
We call a domain $\Omega$ a {\em bounded Lipschitz domain} if 
$\Omega $ is bounded and its boundary $\partial \Omega$
can be covered by a finite number of open balls $B_k$, 
so that, possibly after a proper rotation,
$\partial \Omega \cap B_k$ for each $k$
is a part of the graph of a Lipschitz function. 

\begin{proposition}\label{rychkov}
Let $\Omega \in \Rd$  be a  bounded Lipschitz domain with center $x^0$.
Then there exists a universal bounded linear extension operator  $\ce $ 
for all values of $s,p$, and $q$, i.e., 
\[
(\ce f)|_\Omega = f \qquad \mbox{for all} \quad 
f \in B^s_q (L_p (\Omega))\, , 
\] 
and 
\[
\| \, \ce \, : \, B^s_q (L_p (\Omega)) 
\to B^s_q (L_p (\Rd))\, \| < \infty \, .
\]
In addition we may assume
\begin{equation}\label{eq105}
\supp \, \, \ce f \, \subset \{x \in \Rd : 
\quad |x-x^0| \le 2 \, \diam \, \Omega\}\, .
\end{equation}
\end{proposition}

\begin{rem}
Proposition {\rm \ref{rychkov} }has been proved 
by Rychkov {\rm \cite{Ry}.} Property {\rm (\ref{eq105})} follows from 
Remark {\rm \ref{multiplier}.} 
\end{rem}

Let us now discuss some embedding properties 
of Besov spaces that are needed for our purposes. 

\begin{proposition}\label{compact}
Let $\Omega \subset \Rd$ be an bounded open set. Let 
$0 < p_0,p_1,q_0,q_1 \le \infty$ and let $s,t\in \R$.
Then the embedding
\[
I:  B^{s+t}_{q_0} (L_{p_0} (\Omega)) \to  B^s_{q_1} 
(L_{p_1} (\Omega))
\]
is compact if and only if  
\begin{equation}\label{einb}
t > d \, \Big( \frac{1}{p_0} - \frac{1}{p_1}\Big)_+\, .
\end{equation}
\end{proposition}

\begin{rem}
Sufficiency is proved e.g. in  {\rm \cite{ET96}.} 
The necessity of the given restrictions is almost obvious, 
but see Lemma \ref{einbettung} 
and {\rm \cite{Le99}} for details. 
\end{rem}


Sometimes Besov spaces or Sobolev spaces of 
fractional order are introduced by means of interpolation
(real and/or  complex).
Here we state following, cf. \cite{T02}.
As usual, $(\, \cdot \, , \, \cdot \, )_{\Theta,q}$ 
and $[\, \cdot \, , \, \cdot \, ]_{\Theta}$ 
denote the real and the complex interpolation functor, respectively.

\begin{proposition}\label{interpolation}
Let $\Omega $ be a bounded Lipschitz domain.
Let $0 < q_0,q_1\le \infty$ and let $s_0,s_1 \in \R $.
Let $0 < \Theta < 1$.\\
{\rm (i)} Let $0 < p,q \le \infty$. 
Suppose $s_0 \neq s_1$ and put $s= (1-\Theta)\, s_0 + \Theta\, s_1$.
Then
\[
\Big(B^{s_0}_{q_0} (L_{p}(\Omega)), 
B^{s_1}_{q_1} (L_{p}(\Omega))  \Big)_{\Theta,q} =
B^{s}_{q} (L_{p}(\Omega)) \qquad \mbox{(equivalent quasi-norms)} \, .
\]
{\rm (ii)} Let $0 < p_0,p_1 \le \infty$. 
We put $s= (1-\Theta)\, s_0 + \Theta\, s_1$, 
\[
\frac 1p = \frac{1-\Theta}{p_0} + \frac{\Theta}{p_1} 
\qquad \mbox{and}\qquad 
\frac 1q = \frac{1-\Theta}{q_0} + \frac{\Theta}{q_1} \, .
\]
Then
\[
\Big[ B^{s_0}_{q_0} (L_{p_0}(\Omega)), 
B^{s_1}_{q_1} (L_{p_1}(\Omega))  \Big]_{\Theta} =
B^{s}_{q} (L_{p}(\Omega)) \qquad \mbox{(equivalent quasi-norms)} \, .
\]
\end{proposition}


\subsection{Besov Spaces on Domains -- Intrinsic Descriptions}\label{app5}


For $M\in \Nb$, $h\in \Rd$, and $f: \Rd \to \Cd$ we
define
\[
\Delta_h^M f (x) := \left\{
\begin{array}{lll}
\sum\limits_{j=0}^M {M \choose j} (-1)^{M-j}\, f(x+jh) & 
\qquad & \mbox{if}\quad 
x, \, x+h, \, \ldots \, , x+Mh \in \Omega\, , \\
0 && \mbox{otherwise}\, .
\end{array}
\right.
\]
The corresponding modulus of smoothness is then given by
\[
\omega^M (t,f)_p := \sup_{|h|<t} \, \| \, 
\Delta_h^M f \,  \|_{L_p (\Omega)}\, , \qquad t>0 \, .
\]
The approach by differences coincides with that using 
restrictions as can be seen by the recent 
result of Dispa \cite{Di}. 

\begin{proposition}\label{dispa}
Let $\Omega $ be a bounded Lipschitz domain.  Let $M \in \Nb$.
Let $0< p,q \le \infty$ and  $ d \, \max (0, 1/p-1) < s < M$.
Then 
\begin{eqnarray*}
B^s_q(L_p (\Omega)) & = &  \bigg\{
f \in L_{\max (p,1)} (\Omega): \quad \\
&& 
\| f \|^\Box   :=  \| f\|_{L_p (\Omega)} + \bigg(
\int_0^1 \Big[t^{-s} \, \omega^M (t,f)_p \Big]^q 
\frac{dt}{t} \bigg)^{1/q} <\infty
\bigg\}
\end{eqnarray*}
in the sense of equivalent quasi-norms.
\end{proposition}


\subsection{Sobolev Spaces on Domains}\label{app6}


Let $\Omega $ be a bounded Lipschitz domain.
Let $m \in \Nb$.
As usual $H^m (\Omega)$ denotes the 
collection of all functions $f$ such that the 
distributional derivatives $D^\alpha f$ of order  $|\alpha|\le m$
belong to $L_2 (\Omega)$. The norm is defined as
\[
\| \, f \, |H^m (\Omega)\| 
:= \biggl(  \sum_{|\alpha|\le m} \| \, D^\alpha f\, |L_2(\Omega)\|^2  
\biggr)^{1/2} .
\]
It is well-known that $H^m (\Rd) = B^m_2(L_2(\Rd))$ 
in the sense of equivalent norms,
cf. e.g. \cite{T83}.
As a consequence of the existence of a bounded linear extension operator 
for Sobolev spaces on bounded Lipschitz 
domains, cf. \cite[p.~181]{St}, it follows that 
\[
H^m (\Omega) = B^m_2(L_2(\Omega)) \qquad \mbox{(equivalent norms)}
\]
for such domains. For fractional $s>0$ we introduce the classes by 
complex interpolation. 
Let $0 < s < m$, $s \not\in \Nb$. Then, following \cite[9.1]{LM}, we define
\[
H^s (\Omega) := \Big[H^m (\Omega), L_2(\Omega) \Big]_\Theta\, , 
\qquad \Theta = 1-\frac sm\, .
\]
This definition does not depend on $m$ in the sense of equivalent norms.
This follows immediately from
\[
\Big[ H^{m}(\Omega), L_{2} (\Omega)  \Big]_{\Theta}  = 
\Big[ B^{m}_{2} (L_{2}(\Omega)), B^{0}_{2} (L_{2}(\Omega))  \Big]_{\Theta} =
B^{s}_{2} (L_{2}(\Omega)) \, , \qquad \Theta = 1-\frac sm\, .
\]
(all in the sense of equivalent norms), cf. Proposition \ref{interpolation}.


\subsection{Function Spaces on Domains and Boundary Conditions}
\label{app7}


We concentrate on homogeneous boundary conditions.
Here it makes sense  to introduce two further 
scales of function spaces (distribution spaces).

\begin{definition}
Let $\Omega \subset \Rd$ be an open nontrivial set. 
Let $s \in \R$ and $0 < p,q \le \infty$. 
\\
{\rm (i)}
Then 
$\mathring{B}^s_q(L_p (\Omega))$ denotes the 
closure of ${\mathcal D}(\Omega)$
in $ B^s_q(L_p (\Omega))$, 
equipped with the quasi-norm of  $ B^s_q(L_p (\Omega))$.
\\
{\rm (ii)} Let $s\ge 0$.
Then 
${H}^s_0 (\Omega)$ denotes the closure of ${\mathcal D}(\Omega)$
in $ H^s(\Omega)$, equipped with the norm of  $ H^s(\Omega)$.
\\
{\rm (iii)} By 
$\widetilde{B}^s_{q}(L_p (\Omega))$ we denote the collection of all
$f \in  {\mathcal D}' (\Omega)$ such that 
there is a $g \in B^s_q(L_p (\Rd))$ with
\begin{equation}\label{rand}
g_{\big|{\Omega}} = f \qquad \mbox{and}
\qquad \supp g \subset \overline{\Omega}\, ,
\end{equation}
equipped with the quasi-norm 
\[
\| \, f \, |   
\widetilde{B}^s_q(L_p (\Omega))\| = \inf \| \, g \, | B^s_q(L_p (\Rd))\|\, ,
\]
where the infimum is taken over all such distributions 
$g$ as in {\rm (\ref{rand}).}
\end{definition}

\begin{rem}\label{noch eins}
For a bounded Lipschitz domain 
$\mathring{B}^s_q(L_p (\Omega)) 
=\widetilde{B}^s_q(L_p (\Omega)) = B^s_q(L_p(\Omega))$ holds
if 
\[
 0 <p,q < \infty \, , \quad \max \Big(\frac 1p - 1, \, 
d \, \Big(\frac 1p - 1\Big)\Big) < s < \frac 1p\, , 
\]
cf. {\rm  \cite[Cor.~1.4.4.5]{Gr3}} and {\rm \cite{T02}}. Hence,
\[
H^s_0 (\Omega) = \mathring{B}^s_2(L_2(\Omega)) 
= \widetilde{B}^s_2(L_2 (\Omega)) = {B}^s_2(L_2 (\Omega)) =
H^s (\Omega)
\]
if $0 \le s <1/2$.
\end{rem}

\noindent
Often    it is more convenient to work with a 
scale  $\overline{B}^s_{q}(L_p (\Omega))$, 
originally introduced in \cite{T02}.

\begin{definition}\label{besov*}
Let $\Omega \subset \Rd$ be an open nontrivial set. 
Let $s \in \R$ and $0 < p,q \le \infty$. 
Then we put
\[
\overline{B}^s_q(L_p (\Omega)) := \left\{\begin{array}{lll}
B^s_q(L_p (\Omega) && \qquad \mbox{if}\quad s<1/p\, , 
\\
\widetilde{B}^s_{q}(L_p (\Omega)) && \qquad  \mbox{if} \quad s\ge 1/p\, .
\end{array}\right.
\]
\end{definition}

\noindent
This scale $\overline{B}^s_q(L_p (\Omega))$  
is well-behaved under interpolation and duality, cf. 
\cite{T02}.

\begin{proposition}\label{duality}
Let $\Omega $ be a bounded Lipschitz domain.
Let $1 <p,p_0,p_1,q, q_0,q_1< \infty$ and let $s_0,s_1 \in \R $.
Let $0 < \Theta < 1$.\\
{\rm (i)} 
Suppose $s_0 \neq s_1$ and put $s= (1-\Theta)\, s_0 + \Theta\, s_1$.
Then
\[
\Big(\overline{B}^{s_0}_{q_0} (L_{p}(\Omega)), 
\overline{B}^{s_1}_{q_1} (L_{p}(\Omega))  \Big)_{\Theta,q} =
\overline{B}^{s}_{q} (L_{p}(\Omega)) 
\qquad \mbox{(equivalent quasi-norms)} \, .
\]
{\rm (ii)} 
We put $s= (1-\Theta)\, s_0 + \Theta\, s_1$, 
\[
\frac 1p = \frac{1-\Theta}{p_0} 
+ \frac{\Theta}{p_1} \qquad \mbox{and}\qquad 
\frac 1q = \frac{1-\Theta}{q_0} + \frac{\Theta}{q_1} \, .
\]
Then
\[
\Big[ \overline{B}^{s_0}_{q_0} (L_{p_0}(\Omega)), 
\overline{B}^{s_1}_{q_1} (L_{p_1}(\Omega))  \Big]_{\Theta} =
\overline{B}^{s}_{q} (L_{p}(\Omega)) 
\qquad \mbox{(equivalent quasi-norms)} \, .
\]
{\rm (iii)} With $s\in \R$ and
\[
 1 = \frac{1}{p} + \frac{1}{p'} \qquad \mbox{and}\qquad 
1 = \frac{1}{q} + \frac{1}{q'} 
\]
we find
\[
\Big(\overline{B}^s_{q}(L_p (\Omega))\Big)' 
= \overline{B}^{-s}_{q'}(L_{p'} (\Omega))\, .
\]
Here the duality must be understood in the framework 
of the dual pairing $({\mathcal D}(\Omega),
{\mathcal D}'(\Omega)$.
\end{proposition}


\subsection{Sobolev Spaces with Negative Smoothness} \label{negative}


\begin{definition}
For $s>0$ we define
\[
H^{-s} (\Omega) :=\left\{ \begin{array}{lll}
\Big({H}^s_0 (\Omega)\Big)' & \qquad & \mbox{if}
\quad  s-\frac 12 \neq \mbox{integer}\, , \\
&& \\
 \Big(\widetilde{B}^s_2 (L_2(\Omega))\Big)' && \mbox{otherwise}\, .
\end{array}\right.
\]
\end{definition}

\begin{rem}\label{tilde}
Let $\Omega \subset \Rd$ be a bounded Lipschitz domain. Then 
\[
{H}^s_0 (\Omega) = \widetilde{B}^s_2 (L_2(\Omega))\, , 
\qquad s>0\, , \quad s-\frac 12 \neq \mbox{integer}\, ,
\]
cf. {\rm \cite[Cor.~1.4.4.5]{Gr3}} and 
Proposition {\rm \ref{dispa}.} From 
Remark {\rm \ref{noch eins}} and 
Proposition {\rm \ref{duality}} we conclude the identity
\begin{equation}
H^{-s} (\Omega) = B^{-s}_2 (L_2 (\Omega))\, , \qquad s>0\, , 
\end{equation}
to be understood in the sense of equivalent norms.
\end{rem}

\bigskip
\beginpicture

\setcoordinatesystem units <0.8cm,0.8cm>
\unitlength0.8cm

\put {\vector (0,1) {15} } [Bl] at 1 0
\put {\vector (1,0) {12} } [Bl] at 0 6

\plot 1 6  11 11 /
\plot 1 1  11 6 /

\put {{\small $\widetilde{B}^s_q (L_p (\Omega)) 
= \mathring{B}^s_q (L_p (\Omega))
= {B}^s_q (L_p (\Omega))$}} at 8 6.7
\put {{\small $\widetilde{B}^s_q (L_p (\Omega)) 
= \mathring{B}^s_q (L_p (\Omega))$}} at 4.5 11
\put {{\small $ \neq {B}^s_q (L_p (\Omega))$}} at 5 10

\put {$s$} at 0.5 15
\put {$1/p$} at 12 5.5
\put {$s= 1/p$} at 12 11
\put {$s= 1- 1/p$} at 7.2 15
\put {$s= 1/p - 1$} at 9.2 4
\put {$\bullet$} at 11 6
\put {$1$} at 11 5.5
\put {$\bullet$} at 6 6
\put {$1/2$} at 6.2 5.5
\put {$\bullet$} at 1 12
\put {$1$} at 0.5 12
\put {$\bullet$} at 1 6
\put {$\bullet$} at 1 1
\put {$0$} at 0.5 6.2
\put {$-1$} at 0.5 1


\setdots <0.1cm>

\plot 1 12 6 15 /

\endpicture
\hfill\\

\begin{rem}
{\rm \cite[4.3.2]{T78}.} 
Let $\Omega$ be a bounded open set with a smooth boundary. Then 
$\mathring{B}^s_q(L_p (\Omega)) =\widetilde{B}^s_q(L_p (\Omega))$ holds if  
\[
1 <p,q < \infty\, , \quad \frac 1p - 1  < s < \infty\, , 
\quad s-\frac 1p \neq \mbox{integer}\, .
\]
\end{rem}


\subsection{Wavelet Characterization of Besov Spaces on Domains}\label{app8}


It is a difficult task to construct wavelet bases 
on domains, see \cite[2.12]{C03} and the 
references given there. 
Under certain conditions on the domain $\Omega$ such  
constructions with properties similar
to (\ref{eq203d}), (\ref{eq204}) are known in the literature, see Remark
\ref{r10} above.\\
Let $\Omega$ be a bounded open set in $\Rd$. 
Let  $p,q$ and $s$ be fixed such that
$s > d \max (0, 1/p -1)$. We suppose that
there exist sets $\nabla_j \subset \{ 1,2, \ldots \, ,2^d-1\} \times \Zd$,
with 
\begin{equation}\label{eq205}
0 < \inf_{j=-1, 0, \ldots } 2^{-jd}\,  |\nabla_j| \le 
\sup_{j=-1, 0, \ldots } 2^{-jd}\,  |\nabla_j| < \infty \, ,
\end{equation}
and functions $\psi_{j,\lambda}, \, \widetilde{\psi}_{j,\lambda}  $, 
$\lambda \in \nabla_j$, $j=-1,0,1 , \ldots $, such that
\begin{equation}\label{eq206}
\supp \psi_{j,\lambda}, \quad \supp 
\widetilde{\psi}_{j,\lambda} \subset \Omega \, , \quad 
\lambda \in \nabla_j\, , 
\end{equation}
\begin{equation}\label{eq2051}
\langle \widetilde{\psi}_{i,j,k}, 
\psi_{u,v,\ell} \rangle    =  \delta_{i,u}\,
 \delta_{j,v}\, \delta_{k,\ell}\, , 
\end{equation}  
and such that $f \in B^s_q (L_p (\Omega)) $ if and only if 
\begin{equation}\label{eq207}
f = \sum_{j=-1}^\infty \sum_{\lambda \in \nabla_j}  
\langle f, \widetilde{\psi}_{j,\lambda} 
\rangle \, \psi_{j,\lambda} \qquad 
\mbox{(convergence in ${\mathcal D}'$)}\, , 
\end{equation}
and 
\begin{equation}\label{eq209}
\| \, f \, \|_{B^s_q (L_p(\Omega))}^\clubsuit 
\asymp \| \, f \, \|_{B^s_q (L_p(\Omega))} \, .
\end{equation}
where
\begin{equation}\label{208}
\| \, f \, \|_{B^s_q (L_p(\Omega))}^\clubsuit := \left( \sum_{j=-1}^\infty  
2^{j(s + d(\frac 12 - \frac 1p))q}  \bigg(
\sum_{\lambda \in \nabla_j} 
| \langle f, \widetilde{\psi}_{j,\lambda} \rangle |^p \bigg)^{q/p}
\right)^{1/q} <\infty\, .
\end{equation}


\medskip
\noindent
{\bf Acknowledgment. } \
We thank Stefan Heinrich, Peter Math\'e, Volodya Temlyakov, Hans Triebel,
and Art Werschulz 
for many  valuable remarks and comments. 


\bigskip
\vbox{\noindent Stephan Dahlke\\
Philipps-Universit\"at Marburg\\
FB12 Mathematik und Informatik\\
Hans-Meerwein Stra\ss e\\
Lahnberge\\
35032 Marburg\\
Germany\\
e--mail: {\tt dahlke@mathematik.uni-marburg.de}\\
WWW: {\tt http://www.mathematik.uni-marburg.de/$\sim$dahlke/}\\}

\bigskip
\smallskip
\vbox{\noindent Erich Novak, Winfried Sickel\\
Friedrich-Schiller-Universit\"at Jena\\ 
Mathematisches Institut\\
Ernst-Abbe-Platz 2\\ 
07743 Jena \\ 
Germany\\
e-mail: {\tt \{novak, sickel\}@math.uni-jena.de}\\
WWW: {\tt http://www.minet.uni-jena.de/$\sim$\{novak,sickel\}}}

\bigskip

\end{document}